\documentclass{article}

\usepackage{arxiv}

\usepackage[utf8]{inputenc} 
\usepackage[T1]{fontenc}    
\usepackage{hyperref}       
\usepackage{url}            
\usepackage{booktabs}       
\usepackage{amsfonts}       
\usepackage{nicefrac}       
\usepackage{microtype}      
\usepackage{lipsum}
\usepackage{latexsym}
\usepackage{graphicx}
\usepackage{multicol,multirow}
\usepackage{amsmath,amssymb,amsfonts}
\usepackage{mathrsfs}
\usepackage{amsthm}
\usepackage{apacite}
\usepackage{rotating}
\usepackage{appendix}
\usepackage[authoryear]{natbib}
\usepackage{ifpdf}
\usepackage[T1]{fontenc}
\usepackage{times}
\usepackage{sourcesanspro}
\usepackage{newtxmath}
\usepackage{textcomp}%
\usepackage{xcolor}%
\usepackage{hyperref}
\usepackage{subcaption}

\title{Model Order Reduction based on Runge-Kutta Neural Network}

\author{
  Qinyu Zhuang\\
  Technology\\
  Siemens AG\\
  Munich, Germany \\
  \texttt{qinyu.zhuang@siemens.com} \\
  \And
  Juan Manuel Lorenzi\\
  Technology\\
  Siemens AG\\
  Munich, Germany \\
  \texttt{juan.lorenzi@siemens.com} \\
  \And
  Hans-Joachim Bungartz\\
  Chair of Scientific Computing\\
  Technical University of Munich\\
  Munich, Germany \\
  \texttt{bungartz@tum.de}\\
  \And
  Dirk Hartmann\\
  Technology\\
  Siemens AG\\
  Munich, Germany \\
  \texttt{dirk.hartmann@siemens.com} \\
}

\begin{document}
\maketitle

\begin{abstract}
Model Order Reduction (MOR) methods enable the generation of real-time-capable digital twins, which can enable various novel value streams in industry. While traditional projection-based methods are robust and accurate for linear problems, incorporating Machine Learning to deal with nonlinearity \cite{ubbiali2017reduced, agudelo2009acceleration, pawar2019deep} becomes a new choice for reducing complex problems.
Such methods usually consist of two steps. The first step is dimension reduction by projection-based method, and the second is the model reconstruction by Neural Network. In this work, we apply some modifications for both steps respectively and investigate how they are impacted by testing with three simulation models. In all cases Proper Orthogonal Decomposition (POD) is used for dimension reduction. For this step, the effects of generating the input snapshot database with constant input parameters is compared with time-dependent input parameters. For the model reconstruction step, two types of neural network architectures are compared: Multilayer Perceptron (MLP) and Runge-Kutta Neural Network (RKNN). The MLP learns the system state directly while RKNN learns the derivative of system state and predicts the new state as a Runge-Kutta integrator \cite{devries1994first}.
\end{abstract}


\section{Introduction}
\label{sec:intro}

Physics-based simulation has been an integral part of product development and design as a cheaper alternative to physical prototyping. As the requirements for better simulation performance, the simulation models are usually complex. On the one hand, large computational resource enables the use of such complex models in the design phase. On the other hand, small form-factor hardware is also readily available, enabling "edge computing", i.e. deploying compact computing devices in factories and smart buildings to enable data analysis. The existence of such hardware opens up the possibility of transferring the physics-based models from design phase into the operation phase. This is a key part of the digital twin vision \cite{rasheed2019digital,hartmann2018model}. A digital twin running next to the device can enable novel industrial solutions such as model-based predictive maintenance or model-based optimization.

However, bringing highly complex simulation models into operation phase presents several challenges, which are not present in the design phase. Firstly, the memory footprint of such models needs to be reduced to fit within the limited memory of edge devices, alongside the rest of the other potential data-analysis processes. Secondly, the models need to run at or faster than real-time on the limited hardware of edge devices. The methods that transform simulation models to comply with such requirements are know under the umbrella term model order reduction (MOR)\cite{antoulas2005approximation, hinze2005pod, willcox2002balanced}.

Most MOR techniques rely on a mapping the high-dimensional state of the full order model (FOM) into a lower dimensional space where the reduced order model (ROM) is to be solved. Within this contribution we look into Proper Orthogonal Decomposition (POD), one of the most widely used for finding such mapping \cite{hinze2005pod, willcox2002balanced}. POD is based on performing Principal Component Analysis (PCA) on model trajectories generated by the original FOM, know as snapshots.

For linear models, knowing the dimension-reduction mapping and the equations that define the FOM is enough to obtain a working ROM through projection. However, this cannot be done for complex models with nonlinearity. In this case, a model-reconstruction step is needed in order to reproduce nonlinear behavior in the reduced space. There are several methods perform this step, such as Discrete Empirical Interpolation \cite{chaturantabut2010nonlinear}, Operator Inference \cite{peherstorfer2016data} and Long-short-term-memory Neural Networks \cite{mohan2018deep}. Besides, Recurrent Neural Networks\cite{kosmatopoulos1995high}, have also been extensively used for this purpose \cite{wang2020recurrent, kani2017dr}. As universal approximator, vanilla Multilayer Perceptron (MLP) can be used for this purpose. Besides MLP, we notice a Neural Network model named Runge-Kutta Neural Network (RKNN) which was proposed by Wang and Lin in \cite{wang1998runge} should also be able to work as a surrogate model. This kind of network specializes in non-intrusively modelling the solution of Ordinary Differential Equation (ODE) or Partial Differential Equation (PDE). The RKNN model can simulate the right hand side (RHS) of ODE directly in its sub-networks and step forward in a similar way as a $4^{th}$-order Runge-Kutta integrator, rather than learning the mapping between two successive system states directly.

The combination of POD for dimension reduction and Neural Network represent a purely-data-driven MOR framework. The main contribution of this work consist in exploring methodological variations to this framework. On the one hand, the effects of different approaches to generate the snapshot data for POD are investigated. On the other hand, the impact of using two different Neural Network architectures is explored.

The paper is organized as follows. The introduction for Proper Orthogonal Decomposition is given in Section \ref{sec:pod}. The principle and architecture of Runge-Kutta Neural Network is described in Section \ref{sec:rknn}, numerical experiment evaluating the proposed MOR framework is given in Section \ref{sec:examples}, and their results are shown in Section \ref{sec:results}. Finally, the conclusions follow in Section \ref{sec:conclusions}. As supplement, some figures are provided in Section \ref{appendix}.

\section{Introduction to Proper Orthogonal Decomposition}
\label{sec:pod}

In the last decades, there have been many efforts to develop different techniques in order to obtain compact low dimensional representations of high dimensional datasets. These representations, in general, encapsulate the most important information while discarding less important components. Some applications include image compression using Principal Component Analysis \cite{du2007hyperspectral}, data visualisation using t-SNE \cite{maaten2008visualizing} and structural description of data using Diffusion Maps \cite{coifman2005geometric}. Here the focus will be on Proper Orthogonal Decomposition, which is a variant of Principal Component Analysis in the field of Model Order Reduction.

POD can find a reduced basis $\boldsymbol{V}$ of arbitrary dimension that optimally represents (in a least square sense) the trajectories of the FOM used as snapshots. This basis can be used to project such snapshots into this low-dimensional space.

\subsection{Problem statement}
\label{sec:prob_statement}

The full-dimensional problem that we are intending to reduce in this work is generalized Ordinary Differential Equation (ODE). This ODE is frequently used as governing equation for many engineering problems. Its equation can be written as:

\begin{equation}
    \Dot{\boldsymbol{y}}(t) = \boldsymbol{f}(\boldsymbol{y}(t);\boldsymbol{\mu}(t))
\label{eq:ode_example}
\end{equation}

where $\boldsymbol{y} \in R^{N}$ is the state vector of the FOM. $N$, the number of variables in state vector will be called the size of FOM. And $\mu \in R^{n_\mu}$ is the vector of system parameters which can influence the solution of ODE, where $n_\mu$ is the number of influencing parameters. The target of MOR is to find a reduced model, whose size is $N_r\ll N$, that can reproduce the solution of the FOM with given system parameters to a certain accuracy. This goal will be achieved by mapping the FOM into a reduced space with the help of a reduced basis $\boldsymbol{V}$. As briefly described before, the reduced basis is constructed from the snapshots of the FOM. Snapshots are nothing but a set of solutions, $\boldsymbol{y_1}, \boldsymbol{y_2}, ...,$ $ \boldsymbol{y_{N_s}}$, to the FOM Equation \ref{eq:ode_example} with corresponding system parameter configurations $\boldsymbol{\mu_1}, \boldsymbol{\mu_2}, ..., \boldsymbol{\mu_{Ns}}$. Here, $N_s$ is number of snapshots taken for the FOM. Often, these snapshots will present in the form of matrix called snapshot matrix $\boldsymbol{Y} = [\boldsymbol{y_1}, \boldsymbol{y_2}, ..., \boldsymbol{y_{N_s}}] \in R^{N\times (N_s\cdot k)}$, where $k$ is the number of time steps in each simulation.

\subsection{Taking Snapshots: static-parameter sampling (SPS) and dynamic-parameter sampling (DPS)}
\label{sec:snapshots}

Since the reduced basis $\boldsymbol{V}$ is constructed based on the snapshot matrix $\boldsymbol{Y}$ and later the snapshot will also get involved into training ANN, the quality of the snapshot itself is crucial to the performance of the whole framework. The conventional way of capturing snapshots is designed for data-based intrusive MOR techniques which exploit all information in the internal space of numerical solvers. However, as training data for non-intrusive MOR methods, it shows its limitation of providing more information for model reconstruction. This drives the need of improving snapshot capturing strategy for ANN-based MOR techniques. Here in this paper, we propose a new sampling strategy, \textit{Dynamic Parameter Sampling} (DPS), for taking snapshots to capture the dynamics of the FOM as much as possible.

The simplest way to choose the values for the system parameters $\boldsymbol{\mu}$ is to choose \emph{constant} values $\boldsymbol{\mu(t)} = \boldsymbol{\mu(t_0)}$ using some sparse sampling technique such as \textit{Sparse Grids}, \textit{Latin Hypercube Sampling} \cite{helton2003latin}, \textit{Hammersley Sampling} and \textit{Halton Sampling} \cite{wong1997sampling}. Despite these approaches can often catch enough dynamic modes of the system, the obtained dataset are lacking of diversity of parameter configurations and cannot account for cases in which the parameters are \emph{time-dependent}. Here, we will additionaly explore the use of time-dependent parameters for generating the snapshots.

\begin{equation}
    \Dot{\boldsymbol{y}}_i = \boldsymbol{f}(\boldsymbol{y}_i;\boldsymbol{\mu_i} (t))
\label{eq:td_snapshot}
\end{equation}

In Equation \ref{eq:td_snapshot}, the ODE solved for snapshot $i$ has time-dependent parameter configuration $\boldsymbol{\mu_i}(t)$. This dependency can be any function of time, here to take advantage of sparse sampling techniques, a sinusoidal function is used for dynamic sampling. The following steps explain how DPS works:
\begin{enumerate}
    \item Define the parameter space $\boldsymbol{\mu} \in [\mu_1^{min}, \mu_1^{max}]\times [\mu_2^{min}, \mu_2^{max}] \times ...\times[\mu_{n_\mu}^{min}, \mu_{n_\mu}^{max}]$ for sampling
    \item Use any multivariate sampling to select $N_s$ parameter configurations $\{\boldsymbol{\mu_1, \mu_2, ..., \mu_{N_s}}\}$, where $\boldsymbol{\mu_i}=[\mu_1^i, \mu_2^i, ..., \mu_{n_\mu}^i]$
    \item Derive $N_s$ sinusoidal parameter configurations from the result in last step: $\boldsymbol{\mu_i}(t) = (\boldsymbol{\mu_i}-\boldsymbol{\mu_{min}}) |sin(\omega t)| + \boldsymbol{\mu_{min}}$, where $\boldsymbol{\mu_{min}} = [\mu_1^{min}, \mu_2^{min},  ..., \mu_{n_\mu}^{min}]$ 
    \item Solve the ODE \ref{eq:td_snapshot} with $N_s$ selected parameter configurations and get $N_s$ solutions $\boldsymbol{y_1}$, $\boldsymbol{y_2}$, ..., $\boldsymbol{y_{N_s}}$. Each solution has $k$ time points on the solution trajectory.
\end{enumerate}

\begin{figure}[!htbp]
  \centering
  \label{fig:sta_sampling}{\includegraphics[width=0.45\textwidth,scale=0.5]{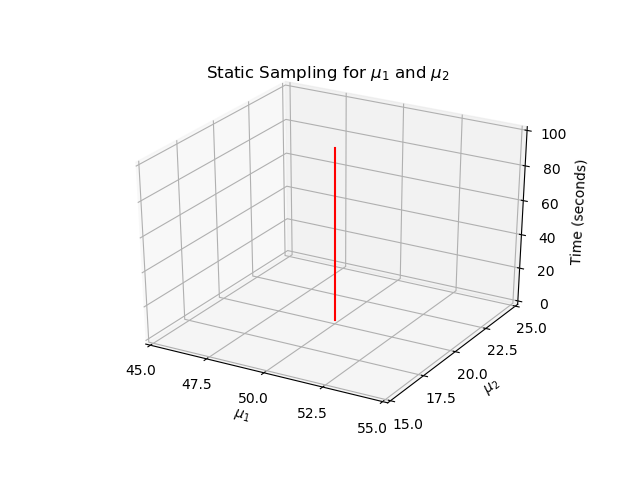}}%
  \hfill
  \label{fig:dyn_sampling}{\includegraphics[width=0.45\textwidth,scale=0.5]{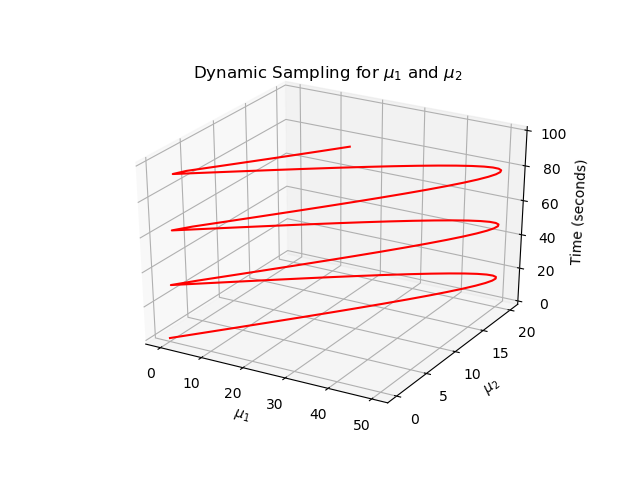}}%
  \hfill
  \captionsetup{font=scriptsize}
  \caption{An example showing the difference between SPS and DPS. We assume parameter space is designed as $[0, 100]\times[0, 100]$. (a) SPS: the parameter configuration $\boldsymbol{\mu} = [50, 20]$ is selected for taking one snapshot. (b) DPS: the parameter configuration $\boldsymbol{\mu} = [50|cos(0.1 t)|, 20|cos(0.1 t)|]$ is selected for running one snapshot simulation}
\label{fig:pod_sampling}
\end{figure}

These $N_s$ solutions are so called snapshots. The snapshots contain the dynamic response of the system with certain parameter configurations. And more dynamic modes are included in the snapshots, more essential information about the system is available to the reduction. Moreover, since the input-output response of the system will be fed into the neural network in Section \ref{sec:rknn}, diversity of the snapshots will strongly influence the training procedure. In this aspect, DPS can provide us with more diverse observation to the system's response.

\subsection{Singular Value Decomposition}
\label{sec:svd}
With the snapshots obtained in the Section \ref{sec:snapshots}, the goal is to find an appropriate reduced basis $\{ \boldsymbol{v_i} \}_{i=1}^{N_r}$ which minimizes the approximation error \cite{chaturantabut2010nonlinear}:

\begin{equation}
\epsilon_{approx} = \sum_{j=1}^{N_s}||\boldsymbol{y}_j - \sum_{i+1}^{N_r}(\boldsymbol{y}_j^T \boldsymbol{v}_i)\boldsymbol{v}_i||_2^2
\label{eq:pod_approximation}
\end{equation}
where $\boldsymbol{y_i}$ stands for $\boldsymbol{y(\mu_i)}$.

The minimization of $\epsilon_{approx}$ in Equation \ref{eq:pod_approximation} can be solved by \textit{Singular Value Decomposition} (SVD) of the snapshot matrix $\boldsymbol{Y} = [\boldsymbol{y_1}, \boldsymbol{y_2}, ..., \boldsymbol{y_{N_s}}] \in R^{N\times (N_s\cdot k)}$:

\begin{equation}
    \boldsymbol{Y} = \boldsymbol{V} \boldsymbol{\Sigma} \boldsymbol{W}^T
\label{eq:svd}
\end{equation}

If the rank of the matrix $\boldsymbol{Y}$ is $k$, then the matrix $\boldsymbol{V}\in R^{N\times k}$ consists of $r$ column vectors $\{\boldsymbol{v}_i, i=1,2,...,k\}$, matrix $\boldsymbol{\Sigma}$ is a diagonal matrix $diag(\sigma_1, \sigma_2, ..., \sigma_k)$. $\sigma_i$ is called $i^{th}$ singular value corresponding to $i^{th}$ singular vector $\boldsymbol{v_i}$.

Essentially, each singular vector $\boldsymbol{v}_i$ represents a dynamic mode of the system. And the corresponding singular value $\sigma_i$ of singular vector $\boldsymbol{v}_i$ can be seen as the "weight" of the dynamic mode $\boldsymbol{v}_i$. Since the greater the weight, the more important the dynamic mode is, $N_r$ singular vectors corresponding to the greatest $N_r$ singular values will be used to construct the reduced basis $\boldsymbol{V}_r = [\boldsymbol{v}_1, \boldsymbol{v}_2, ..., \boldsymbol{v}_{N_r}] \in R^{N\times N_r}$.

Theoretically, the reduced model will have better quality if more dynamic modes are retained in $\boldsymbol{V}$, but this will also increase the size of the reduced model. And since some high frequency noise with small singular value might also be observed in snapshots and later included in the reduced model, the reduced dimension must be selected carefully. There are some research \cite{josse2012selecting, valle1999selection} regarding this topic, and we will not further study it in this work.

\section{Runge-Kutta Neural Network}
\label{sec:rknn}

After finding the projection into the lower dimensional space, we need a way to reproduce the projected dynamics who is governed by Equation \ref{eq:ode_example} in the full order space. In this work, we look into Neural Networks with two architectures named Multilayer Perceptron (MLP) and Runge-Kutta Neural Network (RKNN) respectively.

In both cases the training data uses the snapshots projected into the reduced space.
\begin{align}
\begin{split}
  \mathbf{Y}_r &= \boldsymbol{V^T} \boldsymbol{Y}\\
               &= \boldsymbol{V^T} \left[ \mathbf{y}_{1}, \mathbf{y}_{2}, ..., \mathbf{y}_{N_s} \right]\\
               &= \left[ \mathbf{y}_{r,1}, \mathbf{y}_{r,2}, ..., \mathbf{y}_{r,N_s} \right]
\label{eq:reduced_snapshot_matrix}
\end{split}
\end{align}
where
\begin{equation}
  \mathbf{y}_{r,i} = \left[ \mathbf{y}_{r,i} (t_1), ..., \mathbf{y}_{r,i} (t_k) \right]^T
  \label{eq:reduced_trajectory_matrix}
\end{equation}
corresponds to the projected trajectories corresponding to the parameter $\mu_i(t)$ and $k$ is the number of timesteps used for each of the simulation and if we denote the step size with $\tau$, we will have
\begin{equation}
  t_i = t_0 + \tau j \text{ for } 0 \le j \le k
  \label{eq:time_step_values}
\end{equation}
We also note that $\mathbf{Y}_{r} \in \mathbb{R}^{N_r\times (N_s \cdot k)}$.

\begin{figure}[!htbp]
  \centering
  \includegraphics[scale=0.35]{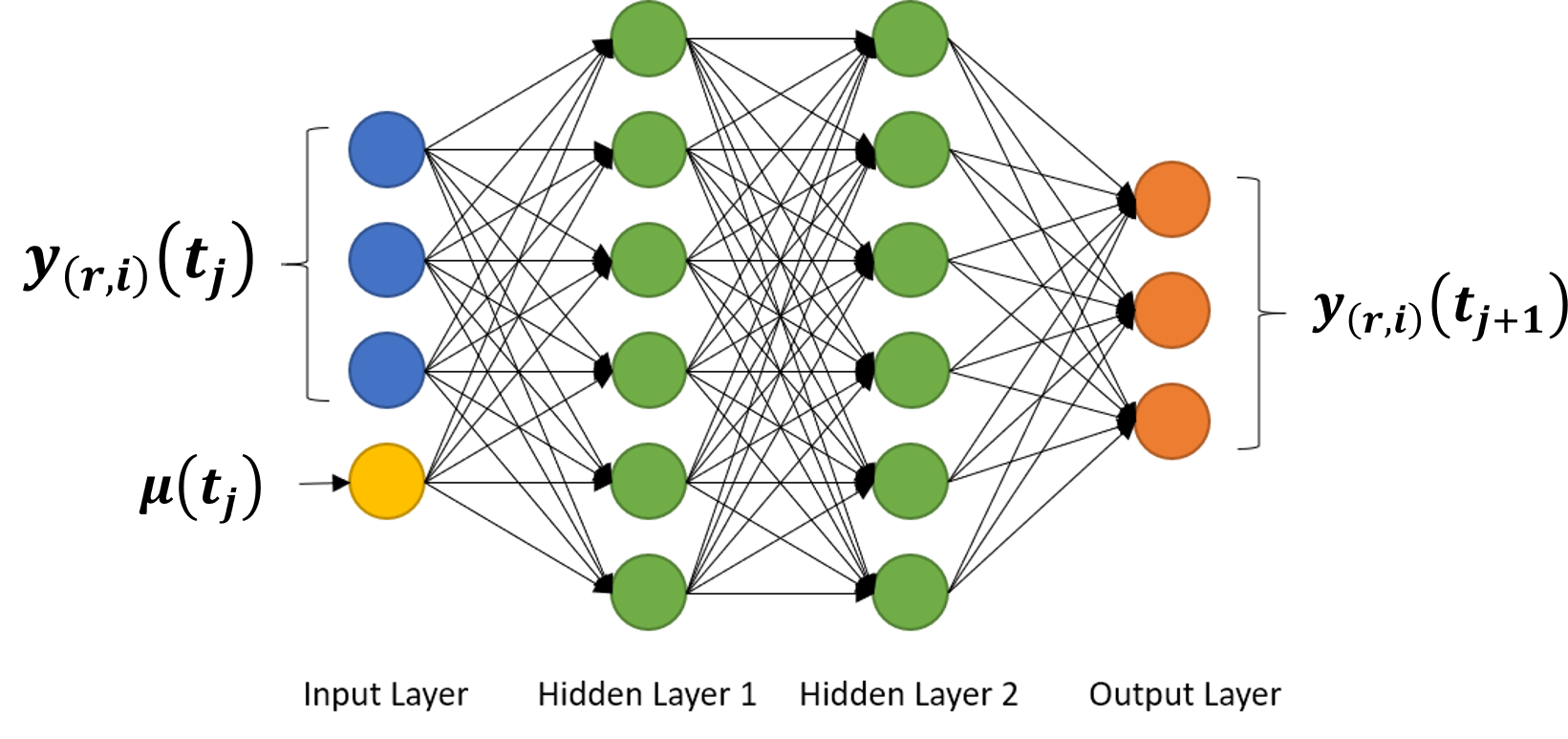}
  \captionsetup{font=scriptsize}
  \caption{Picture of the structure of MLP. MLP predicts the new state $\boldsymbol{y}_{r,i}(t_{j+1})$ based on the given input $\boldsymbol{y}_{r,i}(t_{j})$ and $\boldsymbol{\mu}(t_{j})$}
  \label{fig:mlp}
\end{figure}

An MLP Network learns to predict the evolution of the system by learning the relation between each pair of neighboring state vectors directly as shown in Figure \ref{fig:mlp}.
\begin{equation}
  \mathbf{y}_{r,i} (t_{j+1}) \approx \boldsymbol{g}_{\text{dMLP}} (\mathbf{y}_{r,i} (t_{j}), \boldsymbol{\mu}(t_{j})) \text{ for } 0 \le i \le N_s \text{, } 0 \le j < k
  \label{eq:direct_mlp}
\end{equation}

As alternative,  one can embed a numerical integration scheme into the Neural Network based approach. In this work we focus on the Runge-Kutta 4th order scheme and construct so called RKNN based on it. The $4^th$ order Runge-Kutta integration can be represented as follows:
\begin{equation}
    \boldsymbol{y}_{r,i}(t_{j+1}) = \boldsymbol{y}_{r,i}(t_j) + \frac{1}{6}(\boldsymbol{h}_1 + 2\boldsymbol{h}_2 + 2\boldsymbol{h}_3 + \boldsymbol{h}_4)
\label{eq:rk4th}
\end{equation}
where:
\begin{align*}
    &\boldsymbol{h}_1 = \tau  \boldsymbol{f}(\boldsymbol{y}_{r,i}(t_j); \boldsymbol{\mu}(t_j)), &&\boldsymbol{h}_2 = \tau\boldsymbol{f}(\boldsymbol{y}_{r,i}(t_j) + \frac{\boldsymbol{h}_1}{2};\boldsymbol{\mu}(t_j))\\
    &\boldsymbol{h}_3 = \tau \boldsymbol{f}(\boldsymbol{y}_{r,i}(t_j) + \frac{\boldsymbol{h}_2}{2};\boldsymbol{\mu}(t_j)), &&\boldsymbol{h}_4 = \tau\boldsymbol{f}(\boldsymbol{y}_{r,i}(t_j)+\boldsymbol{h}_3;\boldsymbol{\mu}(t_j)).
\end{align*}
As we can see, the scheme requires the ability to evaluate the right hand side (R.H.S.) $f(\mathbf{y}, \mu)$ four times per time step. The alternative learning scheme is based on using a multilayer perceptron neural network to approximate this function
\begin{equation}
  \boldsymbol{f}(\mathbf{y}; \boldsymbol{\mu}) \approx \boldsymbol{g}_{\text{rkMLP}} (\mathbf{y}; \boldsymbol{\mu})
  \label{eq:approx_rhs}
\end{equation}

Then this MLP can be integrated into the Runge-Kutta scheme as sub-network to obtain an integrator

\begin{equation}
  \boldsymbol{y}_{r,i}(t_{j+1}) = g_{\text{RKNN}} (\boldsymbol{y}_{r,i}(t_{j}); \boldsymbol{\mu}(t_{j})) = \boldsymbol{y}_{r,i}(t_{j}) + \frac{1}{6}(\boldsymbol{h}^{\text{RKNN}}_1 + 2\boldsymbol{h}^{\text{RKNN}}_2 + 2\boldsymbol{h}^{\text{RKNN}}_3 + \boldsymbol{h}^{\text{RKNN}}_4)
\label{eq:rknn}
\end{equation}
where
\begin{align*}
  &\boldsymbol{h}^{\text{RKNN}}_1 = \tau  \boldsymbol{g}_{\text{rkMLP}}(\boldsymbol{y}_{r,i}(t_{j}); \boldsymbol{\mu}(t_{j})), &&\boldsymbol{h}^{\text{RKNN}}_2 = \tau \boldsymbol{g}_{\text{rkMLP}}(\boldsymbol{y}_{r,i}(t_{j}) + \frac{\boldsymbol{h}^{\text{RKNN}}_1}{2};\boldsymbol{\mu}(t_{j}))\\
  &\boldsymbol{h}^{\text{RKNN}}_3 = \tau \boldsymbol{g}_{\text{rkMLP}}(\boldsymbol{y}_{r,i}(t_{j}) + \frac{\boldsymbol{h}^{\text{RKNN}}_2}{2};\boldsymbol{\mu}(t_{j})), &&\boldsymbol{h}^{\text{RKNN}}_4 = \tau \boldsymbol{g}_{\text{rkMLP}}(\boldsymbol{y}_{r,i}(t_{j})+\boldsymbol{h}^{\text{RKNN}}_3;\boldsymbol{\mu}(t_{j})).
\end{align*}

\begin{figure}[!htbp]
  \centering
  \includegraphics[scale=0.5]{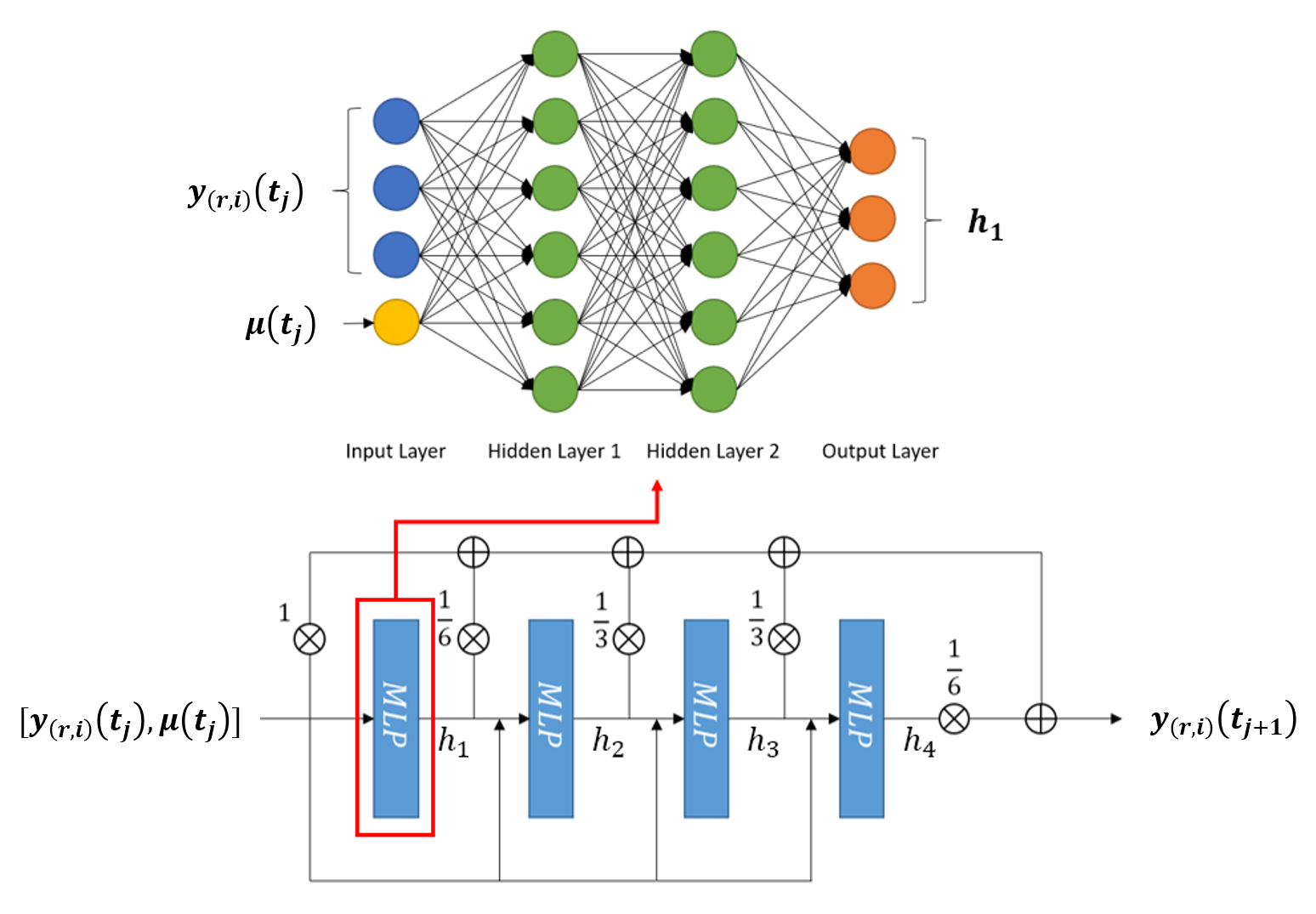}
  \captionsetup{font=scriptsize}
  \caption{Picture of a Runge-Kutta Neural Network (RKNN). RKNN has MLPs as sub-networks which approximate the R.H.S. of Equation \ref{eq:ode_example}. The output of each sub-network is assembled as described in \ref{eq:rk4th}}
  \label{fig:rk_nn}
\end{figure}

In Figure \ref{fig:rk_nn} we present a schematic on how the MLP sub-network is embedded into the larger RKNN based on Equation \ref{eq:rknn}. From Equation \ref{eq:rknn}, we know that an RKNN is essentially a variant of residual network \cite{he2016deep} which learns the increment between two system states instead of learning to map the new state directly from the old state. This leads to a potential advantage of RKNN that we can use deeper network for approximating more complex nonlinearity.

Let's take the computation of $\boldsymbol{h}^{\text{RKNN}}_1$ in Equation \ref{eq:rknn} as an example. We consider an MLP as sub-network with one input layer, 2 hidden layers and one output layer and all fully interconnected. The input layer receives the previous state vector $\boldsymbol{y}_{r,i}(t_{j})$ and the parameter vector $\boldsymbol{\mu_{k}}$, i.e. the number of input neurons is set according to the number of components in the reduced space and the number of parameters. That is, the amount of neurons will depend exclusively on the model and it might change radically for different models. The output from input layer enters two hidden layers using Rectified Linear Unit (ReLU) \cite{li2017convergence} as activation function. Finally, the output from the last hidden layer will be linearly scaled by the output layer and becomes the approximation for $\boldsymbol{h_1}$. The way sub-network computing $\boldsymbol{h}^{\text{RKNN}}_2, \boldsymbol{h}^{\text{RKNN}}_3$ and $\boldsymbol{h}^{\text{RKNN}}_4$ are similar and the only difference is the reduced state vector $\boldsymbol{y}_{r,i}(t_{j})$ in the input will be replaced by $\boldsymbol{y}_{r,i}(t_{j})+c_i\boldsymbol{h}^{\text{RKNN}}_i$, where $c_i$ is the coefficient in Equation \ref{eq:rk4th} and $\boldsymbol{h}^{\text{RKNN}}_i$ is computed in previous iteration. The parameters of the sub-network are updated using backpropagation \cite{rumelhart1986learning}.

For the activation function, there are also multiple choices besides ReLU, e.g. Sigmoid or tanh \cite{karlik2011performance}. The decision should be made for the Neural Network individually. But for most cases ReLU produced faster convergence. Due to this it is suggested to use ReLU as activation function for all hidden layers. A problem of ReLU units is "dead" neuron. That means once the input of a neuron is 0, due to "zero-gradient" in the non-positive region of ReLU function, the parameters of this neuron will never be updated. Some dead neurons can give Neural Network suitable sparsity. But it can be trouble if most of neurons are dead. In this case a variant of ReLU called Leaky ReLU \cite{dubey2019comparative} can be considered as alternative.

Similar to explicit integrators, implicit ones can be used \cite{martinez1993continuous}. Due to their recurrent nature, or the dependence of predictions on themselves, the architecture of this type of neural network will have recurrent connections and offer an alternative to the approach taken here.

\section{Numerical Examples}
\label{sec:examples}
In this section, the three numerical examples used for testing are described. The first one is a computer heat sink, which is governed by thermal equation without any nonlinear contribution. The second one is a gap-radiation model whose physics equation is thermal equation including $4^{th}$-order nonlinear radiation terms. And the last example is a thermal-fluid model simulating a heat exchanger, which includes the effects of fluid flow.

All test cases are all firstly reduced by POD method and then ROMs are evaluated by both MLP and Runge-Kutta neural networks. And to validate the ROMs, each example will be assigned two test cases, one is a constant-load test dynamic-load test.

The numerical experiments are performed on a computer with Intel Xeon E5-2640 CPU and 48GB Memory. The reference solution is provided by Thermal/Flow Multi-physics module in simulation software Siemens Simcenter NX 12 \footnote{https://www.plm.automation.siemens.com/global/en/products/simcenter/}.

\subsection{Heat Sink Model}
The first FEM model simulates a chip cooled by attached heat sink. The chip will be the heat source of the system and heat flux travels from the chip to the heat sink then released into the environment through fins. Since all material used in the model is not temperature-dependent, the system is governed by the linear heat transfer equation as Equation \ref{eq:heatsink_eq}.

\begin{equation}
    \boldsymbol{C_p} \Dot{\boldsymbol{T}} = \boldsymbol{K T} + \boldsymbol{B u}
    \label{eq:heatsink_eq}
\end{equation}
where $\boldsymbol{C_p}$ is thermal capacity matrix, $\boldsymbol{K}$ is thermal conductivity matrix, $\boldsymbol{B}$ is generalized load matrix, $\boldsymbol{T}$ is state vector of temperature and $\boldsymbol{u}$ is the input heat load vector.

\begin{figure}[htbp]
  \centering
  \begin{subfigure}[b]{0.45\textwidth}
  \centering
  \includegraphics[width=\textwidth]{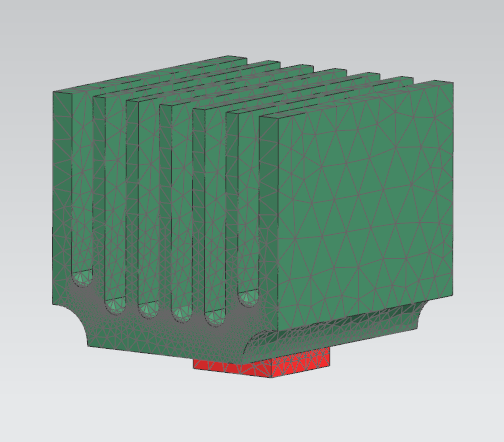}
  \captionsetup{font=scriptsize}
  \caption{Picture of the heat sink model}
  \label{fig:heatsink_model}
  \end{subfigure}
  \hfill
  \begin{subfigure}[b]{0.45\textwidth}
  \centering
  \includegraphics[width=\textwidth]{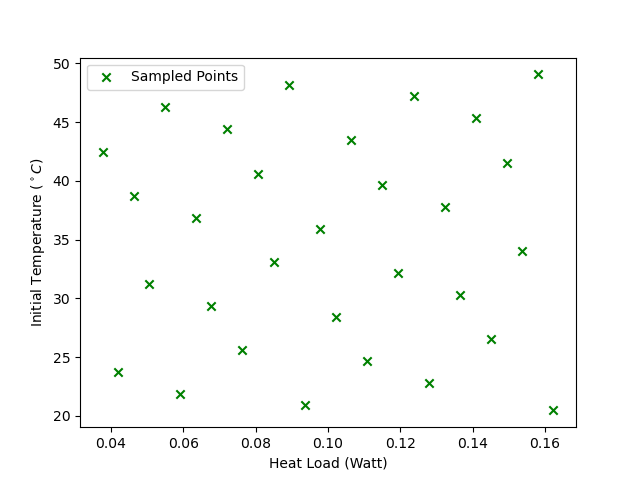}
  \captionsetup{font=scriptsize}
  \caption{Sampled parameters}
  \label{fig:sampled_data_heatsink}
  \end{subfigure}
  \hfill
 \captionsetup{font=footnotesize}
 \caption{Test case: heat sink model}
\end{figure}

The meshed model is sketched in Figure \ref{fig:heatsink_model}. The model consists of two parts, the red part is a chip made from copper whose thermal capacity is 385 $J/(kg\cdot K)$ and thermal conductivity is 387 $W/(m\cdot K)$. The green part represents the attached heat sink made from aluminum. The relevant material properties are 963 $J/(kg\cdot K)$ thermal capacity and 151 $W/(m\cdot K)$ thermal conductivity.

In the designed space there are two variable parameters: initial temperature $T_0$ and input heat load $u$. The range of initial temperature is from 20 $^\circ C$ to 50 $^\circ C$. The power of the volumetric heat source in the chip has range from 0.05 $W/mm^3$ to 0.15 $W/mm^3$. Parameter coordinates sampled in parameter space are shown in Figure \ref{fig:sampled_data_heatsink}.

Both test cases have uniform initial temperature of 20 $^\circ C$. The volumetric heat source used by constant-load test has fixed power 0.1 $W/mm^3$ and the one used by dynamic-load test has power as function of time $u(t) = 0.15 - 0.1 \frac{t}{500}$ $W/mm^3$.

\subsection{Gap-radiation Model}
The model presented here is a thermal model including radiative coupling. This model can be though of as a proxy for other radiation-dominated heat transfer models such a those occurring in aerospspace, energy and manufacturing. The discrete governing equation is as Equation \ref{eq:radiation_eq}.

\begin{equation}
    \boldsymbol{C_p} \Dot{\boldsymbol{T}} = \boldsymbol{K T} + \boldsymbol{R T}^4 + \boldsymbol{B u}
    \label{eq:radiation_eq}
\end{equation}
in addition to Equation \ref{eq:heatsink_eq}, $\boldsymbol{R}$ is the radiation matrix.

\begin{figure}[htbp]
  \centering
  \begin{subfigure}[b]{0.45\textwidth}
  \centering
  \includegraphics[width=\textwidth]{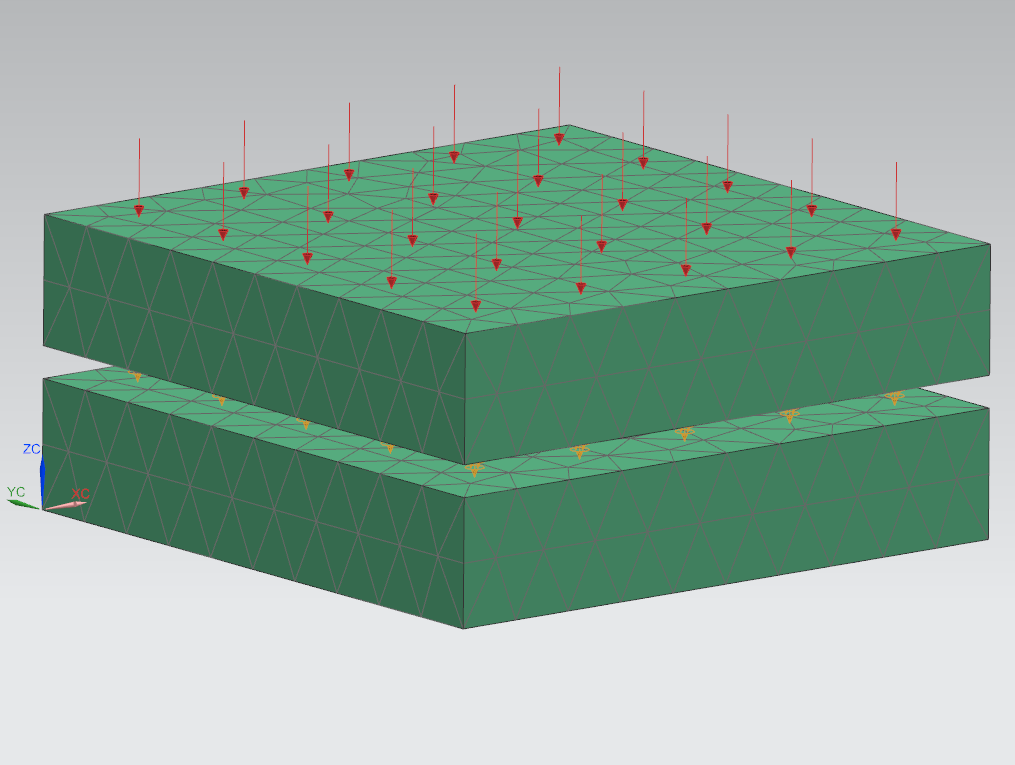}
  \captionsetup{font=scriptsize}
  \caption{Picture of the gap-radiation model}
  \label{fig:radiation_model}
  \end{subfigure}
  \hfill
  \begin{subfigure}[b]{0.45\textwidth}
  \centering
  \includegraphics[width=\textwidth]{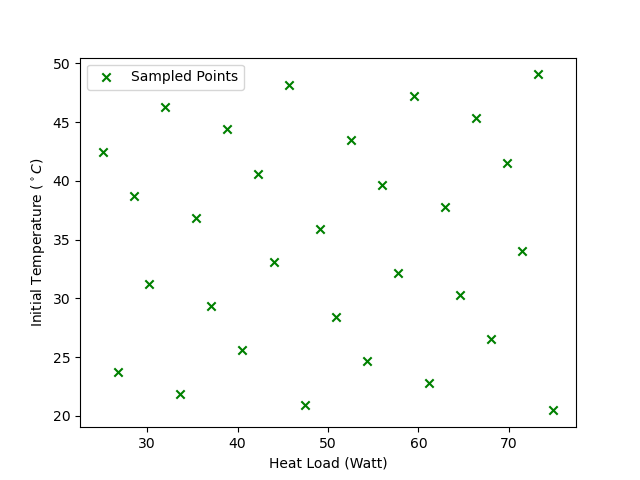}
  \captionsetup{font=scriptsize}
  \caption{Sampled parameters}
  \label{fig:sampled_data_radiation}
  \end{subfigure}
  \hfill
\end{figure}

The 3D model is sketched in Figure \ref{fig:radiation_model}. As shown, The upper plate is heated by thermal flow. When the temperature equilibrium between the two plates is broken due to an increment of the temperature of the upper plate, radiative flux due to the temperature difference between two surface of the gap takes place.

Both plates are made of steel, whose thermal capacity is 434 $J/(kg\cdot K)$ and thermal conductivity is 14 $W/(m\cdot K)$. The initial temperature $T_0$ and the applied heat load $u$ are the variable parameters. The time scale of the whole simulated process is from 0 to 3600 seconds.

Assuming operation condition where the input heat load of the system has the lower limitation $u_{min}=40W$ and upper limitation $u_{max}=60W$ and the initial temperature has the lower limitation $T_{0,min}=20^\circ C$ and upper limitation $T_{0,max}=300^\circ C$ is investigated. So the parameter space to take snapshots is defined as $[20^\circ C, 300^\circ C]\times[40 W, 60 W]$. And parameter configurations used for snapshots are given in Figure \ref{fig:sampled_data_radiation}.

The test scenarios includes a constant-load test and a dynamic-load test. The constant-load test has fixed load magnitude of 50 $W$, while the dynamic-load test has a time-dependent load magnitude $u(t) = 60 - 20 \frac{(t-1800)^2}{1800^2}$ $W$.

\subsection{Heat Exchanger Model}
\label{sec:exchanger}

\begin{figure}[htbp]
  \centering
  \begin{subfigure}[b]{0.45\textwidth}
  \centering
  \includegraphics[width=\textwidth]{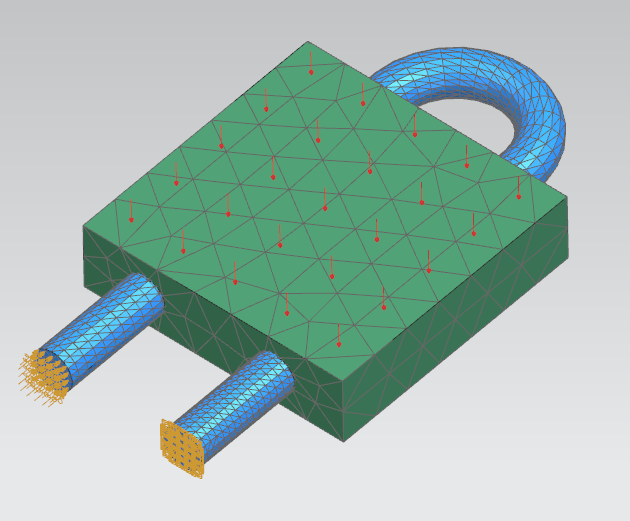}
  \captionsetup{font=scriptsize}
  \caption{Picture of the heat exchanger model}
  \label{fig:exchanger_model}
  \end{subfigure}
  \hfill
  \begin{subfigure}[b]{0.45\textwidth}
  \centering
  \includegraphics[width=\textwidth]{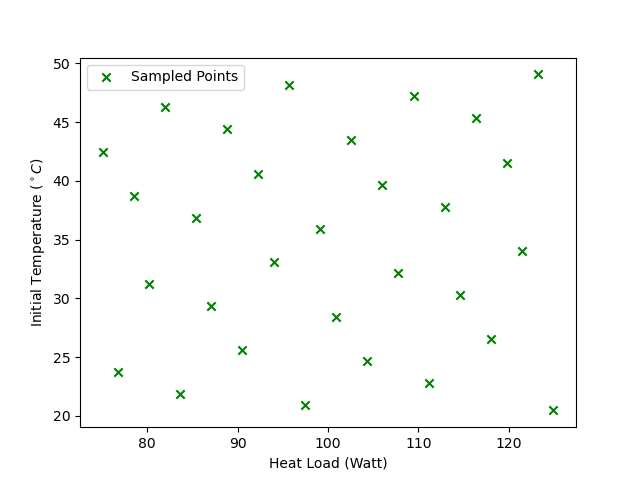}
  \captionsetup{font=scriptsize}
  \caption{Sampled parameters}
  \label{fig:sampled_data_exchanger}
  \end{subfigure}
  \hfill
\end{figure}

In this example, we model the temperature profile of a solid that is in contact with a fluid channels to cool it down. The flow-field is governed by  the Navier-Stokes equations written in integral form \cite{kassab2003bem}:

\begin{equation}
\label{eq:ns_equation}
 \int_{\Omega} \frac{\partial W}{\partial t} d\Omega + \int_{\Gamma} (F-T)\cdot n d\Gamma = \int_{\Omega} Sd\Omega
\end{equation}
where $\Omega$ denotes the volume, $\Gamma$ denotes the surface bounded by the volume $\Omega$, and $n$ is the outward-drawn normal. The conserved variables are contained in $W=(\rho, \rho u, \rho v, \rho w, \rho e, \rho k, \rho \omega)$, where, $\rho, u, v, w, e, k, \omega$ are the density, the velocity in x-, y- and z-directions, and the specific total energy. $F$ and $T$ are convective and diffusive fluxes, respectively, $S$ is a vector containing all terms arising from the use of a non-inertial reference frame as well as in the production and dissipation of turbulent quantities.

The governing equation for heat-conduction field in the solid is:
\begin{equation}
    \nabla\cdot[k_s(T_s)\nabla T_s]=0
    \label{eq:heatconduction_equation}
\end{equation}
here $T_s$ denotes the temperature of the solid, and $k_s$ is the thermal conductivity of the solid material.

Finally the equations should be satisfied on the interface between fluid and solid are:
\begin{equation}
    T_f = T_s \\
    k_f\frac{\partial T_f}{\partial n} = -k_s \frac{\partial T_s}{\partial n} \\
    \label{eq:interface_condition}
\end{equation}
here $T_f$ is the temperature computed from N-S solution Equation \ref{eq:ns_equation} and $k_f$ is thermal conductivity of fluid.

The liquid flowing in the cooling tube is water whose thermal capacity is 4187 $J/(kg\cdot K)$ and thermal conductivity is 0.6 $W/(m\cdot K)$. The fluid field has an inlet boundary condition which equals to $1m/s$. The solid piece is made of aluminum whose thermal capacity is 963 $J/(kg\cdot K)$ and thermal conductivity is 151 $W/(m\cdot K)$.

The initial temperature in the designed parameter space is from 20 $^\circ C$ to 50 $^\circ C$ and the heat load applied on the solid component varies from 80 to 120 $W$, as shown in Figure \ref{fig:sampled_data_exchanger}. Therefore, we will use a constant heat load whose magnitude is 100 $W$ and a variant heat load whose magnitude is a function of time $u(t) = 120 - 0.8 t$ $W$ to validate the ROM respectively.

\section{Results}
\label{sec:results}
In this section, the results of the numerical tests in Section \ref{sec:examples} are presented. The comparison will be made along two dimension: the first dimension is sampling method, i.e. Static Parameter Sampling (SPS) versus Dynamic Parameter Sampling (DPS).
The other dimension is architecture of neural network, i.e. Multilayer Perceptron versus Runge-Kutta Neural Network. Some qualitative and quantitative conclusion will be drawn based on comparison.

In Section \ref{sec:pod_rslt}, we will investigate how sampling strategy influences construction of reduced basis. And in Section \ref{sec:sps_vs_dps} and Section \ref{sec:mlp_vs_rknn}, we will focus on how dataset and architecture influences the learning quality respectively.

\subsection{POD Reduction}
\label{sec:pod_rslt}
With SPS and DPS respectively, 30 snapshots are taken for each system with 30 different heat-load magnitude and initial temperature. Each snapshot is transient analysis from $0$ to $T_{end}$ with 100 time points solved by NX 12.

For DPS, $\omega=\pi/T_{end}$ is chosen to be a sinusoidal parameter function. In this case, parameter configuration $\boldsymbol{\mu}$ will change sinusoidally between $\boldsymbol{\mu}_{min}$ and $\boldsymbol{\mu}_{max}$. Both values are defined by SPS.

\begin{figure}[htbp]
  \centering
  \includegraphics[width=.6\textwidth]{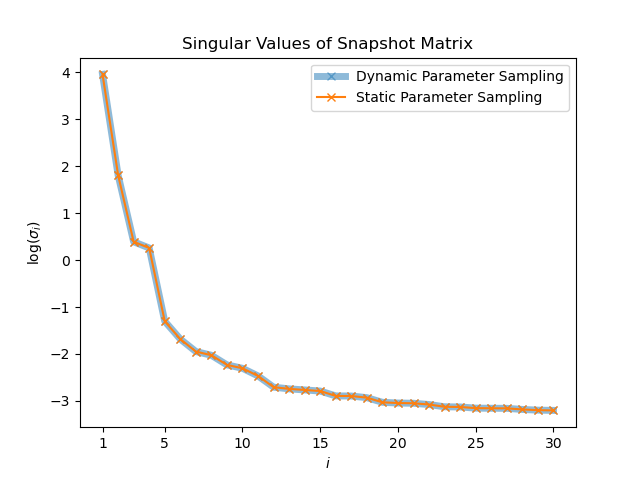}
  \captionsetup{font=scriptsize}
  \caption{Heat sink model: Singular values of SPS snapshots and DPS snapshots}
  \label{fig:sigmas_sink}
\end{figure}

\begin{figure}[htbp]
  \centering
  \includegraphics[width=.6\textwidth]{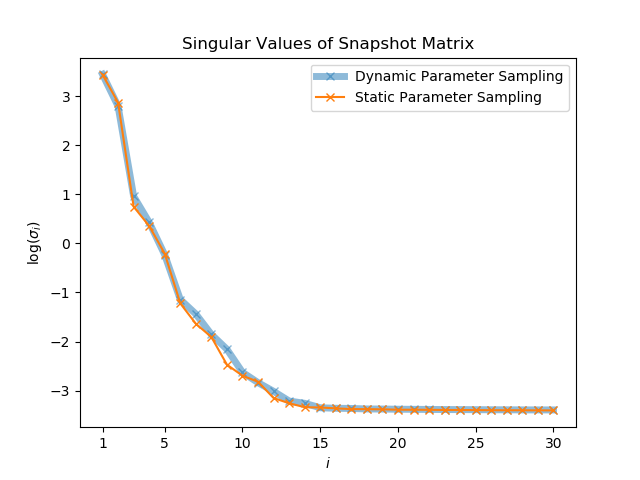}
  \captionsetup{font=scriptsize}
  \caption{Gap radiation model: Singular values of SPS snapshots and DPS snapshots}
  \label{fig:sigmas_gap}
\end{figure}

\begin{figure}[htbp]
  \centering
  \includegraphics[width=.6\textwidth]{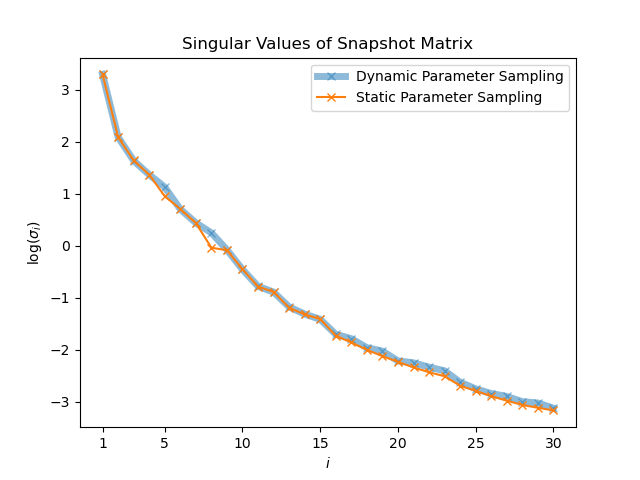}
  \captionsetup{font=scriptsize}
  \caption{Heat exchanger model: Singular values of SPS snapshots and DPS snapshots}
  \label{fig:sigmas_exchanger}
\end{figure}

The first 30 singular values found from SPS snapshots and DPS snapshots are shown in Figure \ref{fig:sigmas_sink}, Figure \ref{fig:sigmas_gap} and Figure \ref{fig:sigmas_exchanger} for each test case respectively. In Heat sink test model, the curves of $\sigma$ from calculated from different sampling strategies seem to exactly overlap on each other. And in Figure \ref{fig:sigmas_gap} and Figure \ref{fig:sigmas_exchanger}, only slight difference can be found between two curves.

However, this does not mean the reduced basis constructed by different dataset has the same quality. To quantitatively measure the quality of reduced basis, here we perform re-projection test to calculate the error generated by each reduced basis. In the test, the reference trajectory is firstly projected into then reduced space then back-projected into the full space. This process can be done by using $\boldsymbol{\hat{Y}}_{ref} = \boldsymbol{V}\boldsymbol{V}^T\boldsymbol{Y}_{ref}$. Depending on re-projected trajectory, we can calculate the relative re-projection error using $E_{pod}=\frac{1}{N_s\cdot k}\frac{|\boldsymbol{Y}_{ref}-\boldsymbol{\hat{Y}}_{ref}|_2}{|\boldsymbol{Y}_{ref}|_2}$, where $N_s$ is number of snapshot simulations and $k$ is the number of time points in each simulation.

The POD error with different size of reduced basis is shown in Figure \ref{fig:pod_error_heatsink} - \ref{fig:pod_error_radiation}. Similar to summarized above, in most situations, POD error curves generated by different sampling strategies can converge to similar values when the size of ROM is large enough. Nevertheless, in Figure \ref{fig:poder_heatsink_varnt}, 
the re-projection error of SPS-ROM converges to a value greater than DPS-ROM's. This means the reduced basis calculated by DPS snapshots has higher quality.

It is observed that, SPS-ROM can have smaller projection error in some cases.This suggest that hybrid datasets (including static and dynamic parameters) may be needed to improve this aspect even further. We propose to study this possibility in future works.

\subsection{ANN-based ROM}
\label{sec:rom_ann}

The quality of the ROM is affected by the reduced basis constructed by POD and the neural network trained in the reduced space. In Section \ref{sec:pod_rslt} by comparing singular values calculated from snapshot matrices sampled in different approaches, we did not see different sampling strategies have huge impact on discovering latent physics. But they can still make difference on training neural network. This will be investigated in Section \ref{sec:sps_vs_dps}. And the performance of RKNN in the reduced space will be shown in Section \ref{sec:mlp_vs_rknn}.

\subsubsection{Static-Parameter Sampling versus Dynamic-Parameter Sampling}
\label{sec:sps_vs_dps}

In Section \ref{sec:snapshots}, we introduced two approaches for taking snapshots. In conventional static-parameter sampling, the parameter configuration assigned to the system is constant during the simulation time, i.e. $\boldsymbol{\mu}(t) = \boldsymbol{\mu}_0$. To include more parameter diversity in training data, we designed new sampling strategy named dynamic-parameter sampling whose parameter configuration is function of time. Datasets sampled by different strategies influences not only the quality of reduced basis but also the training of neural network. In this section, how different datasets can affect the training is investigated.

To measure the prediction quality of neural network correctly, the influence of reduction projection will be removed from measurement. In other words, the approximation error should be evaluated in the reduced space. Therefore, the reference trajectory is reduced by $\boldsymbol{Y_r}_{ref} = \boldsymbol{V^T}\boldsymbol{Y}_{ref}$. With the prediction $\boldsymbol{\Tilde{Y_r}}$ made by neural network, the relative approximation error is calculated as $E_{ann}=\frac{1}{N_s\cdot k}\frac{|\boldsymbol{Y_r}_{ref}-\boldsymbol{\Tilde{Y_r}}|_2}{|\boldsymbol{Y_r}_{ref}|_2}$, where $N_s$ is number of snapshot simulations and $k$ is the number of time points in each simulation.

By observing Table \ref{tab:hs_spsvsdps} - \ref{tab:he_spsvsdps}, in most cases ROMs trained by DPS data have better prediction in both test cases. But unlike many other MOR methods, the simulation error does not monotonically decrease while increasing the size of ROM. When the ROM size is very small, the simulation error is mainly caused by POD projection error. In this region, simulation error will decrease with increasing ROM size. Because the most dominant POD components are the most important physics modes of the system, by adding them the ROM knows more about the original system. However, POD components with smaller singular values might be trivial information in snapshots, e.g. noise from numerical integration. These trivial information also becomes noise in training data and disturbs the learning process. This also reminds us to be careful to choose the size of ROM while using ANN-based MOR technique.

The tests here are considered as independent test since test datasets are provided by new NX simulation with different parameter configuration from the simulations used in the snapshots. The independent test is optimal if the two datasets originate from two different sampling strategies \cite{ozesmi2006methodological}. Here dynamic-load test can be considered to be optimal independent test for SPS-ROM and vice versa for DPS-ROM. As shown in Table \ref{tab:hs_spsvsdps} - \ref{tab:he_spsvsdps}, the ANNs trained by DPS data have better performance on its optimal independent tests. But ANNs trained by SPS data also outperforms in some test cases. This result reminds us again that using dataset obtained from diverse sources should be considered.

\begin{table}[htbp]
\small
 \begin{subtable}{.45\linewidth}
 \centering
 \captionsetup{font=footnotesize}
 \caption{Constant-load Test}
 \begin{tabular}{ccccc}
  \toprule
  $N_r$ & MLP + SPS & MLP + DPS\\
  \midrule
1 & 0.09 & 0.16($+0.07$)\\
2 & 0.05 & 0.19($+0.14$)\\
3 & 0.11 & 0.13($+0.02$)\\
4 & 0.10 & 0.16($+0.06$)\\
5 & 0.30 & 0.18($-0.12$)\\
10 & 0.19 & 0.19($-0.00$)\\
15 & 0.10 & 0.16($+0.06$)\\
20 & 0.11 & 0.15($+0.04$)\\
25 & 0.19 & 0.16($-0.04$)\\
30 & 0.17 & 0.16($-0.01$)\\
  \bottomrule
 \end{tabular}
 \label{tab:hs_const_test}
 \end{subtable}
 \hfill
 \begin{subtable}{.45\linewidth}
 \centering
 \captionsetup{font=footnotesize}
 \caption{Dynamic-load Test}
 \begin{tabular}{ccccc}
  \toprule
  $N_r$ & MLP + SPS & MLP + DPS\\
  \midrule
1 & 0.38 & 0.15($-0.23$)\\
2 & 0.34 & 0.20($-0.15$)\\
3 & 0.33 & 0.15($-0.17$)\\
4 & 0.31 & 0.18($-0.13$)\\
5 & 0.37 & 0.18($-0.19$)\\
10 & 0.34 & 0.19($-0.15$)\\
15 & 0.29 & 0.17($-0.12$)\\
20 & 0.31 & 0.17($-0.14$)\\
25 & 0.32 & 0.16($-0.16$)\\
30 & 0.31 & 0.17($-0.14$)\\
  \bottomrule
 \end{tabular}
 \label{tab:hs_varnt_test}
 \end{subtable}
\caption{SPS versus DPS: heat sink model, Relative error(\%)}
\label{tab:hs_spsvsdps}
\end{table}

\begin{table}[htbp]
\small
 \begin{subtable}{.45\linewidth}
 \centering
 \captionsetup{font=footnotesize}
 \caption{Constant-load Test}
 \begin{tabular}{ccccc}
  \toprule
  $N_r$ & MLP + SPS & MLP + DPS\\
  \midrule
1 & 0.50 & 0.59($+0.09$)\\
2 & 0.59 & 0.27($-0.32$)\\
3 & 0.25 & 0.28($+0.03$)\\
4 & 0.17 & 0.27($+0.10$)\\
5 & 0.22 & 0.18($-0.04$)\\
10 & 0.13 & 0.14($+0.01$)\\
15 & 0.17 & 0.22($+0.05$)\\
20 & 0.22 & 0.13($-0.10$)\\
25 & 0.19 & 0.18($-0.01$)\\
30 & 0.16 & 0.16($-0.00$)\\
  \bottomrule
 \end{tabular}
 \label{tab:gr_const_test}
 \end{subtable}
 \hfill
 \begin{subtable}{.45\linewidth}
 \centering
 \captionsetup{font=footnotesize}
 \caption{Dynamic-load Test}
 \begin{tabular}{ccccc}
  \toprule
  $N_r$ & MLP + SPS & MLP + DPS\\
  \midrule
1 & 0.47 & 0.52($+0.05$)\\
2 & 0.63 & 0.19($-0.44$)\\
3 & 0.30 & 0.26($-0.04$)\\
4 & 0.25 & 0.15($-0.10$)\\
5 & 0.17 & 0.19($+0.02$)\\
10 & 0.27 & 0.19($-0.08$)\\
15 & 0.21 & 0.16($-0.04$)\\
20 & 0.21 & 0.13($-0.08$)\\
25 & 0.17 & 0.13($-0.04$)\\
30 & 0.19 & 0.12($-0.07$)\\
  \bottomrule
 \end{tabular}
 \label{tab:gr_varnt_test}
 \end{subtable}
\caption{SPS versus DPS: gap radiation model, Relative error(\%)}
\label{tab:gr_spsvsdps}
\end{table}

\begin{table}[htbp]
\small
 \begin{subtable}{.45\linewidth}
 \centering
 \captionsetup{font=footnotesize}
 \caption{Constant-load Test}
 \begin{tabular}{ccccc}
  \toprule
  $N_r$ & MLP + SPS & MLP + DPS\\
  \midrule
1 & 0.20 & 0.27($+0.07$)\\
2 & 0.31 & 0.27($-0.04$)\\
3 & 0.25 & 0.25($-0.00$)\\
4 & 0.24 & 0.24($-0.00$)\\
5 & 0.13 & 0.16($-0.03$)\\
10 & 0.17 & 0.11($-0.05$)\\
15 & 0.12 & 0.16($+0.04$)\\
20 & 0.14 & 0.12($-0.02$)\\
25 & 0.15 & 0.15($-0.00$)\\
30 & 0.17 & 0.13($-0.04$)\\
  \bottomrule
 \end{tabular}
 \label{tab:he_const_test}
 \end{subtable}
 \hfill
 \begin{subtable}{.45\linewidth}
 \centering
 \captionsetup{font=footnotesize}
 \caption{Dynamic-load Test}
 \begin{tabular}{ccccc}
  \toprule
  $N_r$ & MLP + SPS & MLP + DPS\\
  \midrule
1 & 0.48 & 0.22($-0.26$)\\
2 & 0.38 & 0.17($-0.21$)\\
3 & 0.28 & 0.27($-0.02$)\\
4 & 0.31 & 0.21($-0.11$)\\
5 & 0.17 & 0.15($-0.02$)\\
10 & 0.22 & 0.13($-0.09$)\\
15 & 0.17 & 0.17($-0.00$)\\
20 & 0.18 & 0.14($-0.04$)\\
25 & 0.15 & 0.15($-0.00$)\\
30 & 0.16 & 0.14($-0.02$)\\
  \bottomrule
 \end{tabular}
 \label{tab:he_varnt_test}
 \end{subtable}
\caption{SPS versus DPS: heat exchanger model, Relative error(\%)}
\label{tab:he_spsvsdps}
\end{table}

\subsubsection{Multilayer Perceptron versus Runge-Kutta Neural Network}
\label{sec:mlp_vs_rknn}

To further improve the capability of predicting ODE systems, we propose to use network architecture inspired by Runge-Kutta integrator in Section \ref{sec:rknn}. Although it has been previously seen that RKNN can be used to make long-term predictions for ODE system \cite{wang1998runge}, it is still not clear if such approach can work for POD-projected ODE systems. A possible challenge arises from the highly variable scale differences between the reduced dimensions. Typically, the reduced snapshot will have much larger component in some of the reduced dimensions (those with large singular values) than in others (small singular values). This presents a challenge for NN-based learning methods.

The structure of MLP and RKNN is the same as described in Figure \ref{fig:mlp} and Figure \ref{fig:rk_nn}, but in detail, the number of input neurons $n_{in}$ depends on the size of ROM and number of system parameters, i.e. $n_{in} = N_r + n_\mu$. The number of neurons on each hidden layers is chosen to be $32$, this value is always greater than input features, i.e. $N_r + n_\mu$. This is important for neural network to interpret the underlying dynamics. And same as depicted in the figures, each network has 1 input layer, 2 hidden layers and 1 output layer, which enables network to approximate any kind of mathematical function.

The prediction error $E_{ann}$ is measured in the same way as in section \ref{sec:sps_vs_dps}. Since it is known from the last comparison that DPS dataset is considered to be better training dataset, here both types of neural networks are trained with DPS dataset.

\begin{table}[htbp]
\small
 \begin{subtable}{.45\linewidth}
 \centering
 \captionsetup{font=footnotesize}
 \caption{Constant-load Test}
 \begin{tabular}{ccccc}
  \toprule
  $N_r$ & MLP + DPS & RKNN + DPS\\
  \midrule
1 & 0.16 & 0.16($-0.01$)\\
2 & 0.19 & 0.15($-0.04$)\\
3 & 0.13 & 0.12($-0.01$)\\
4 & 0.16 & 0.16($-0.01$)\\
5 & 0.18 & 0.17($-0.01$)\\
10 & 0.19 & 0.17($-0.02$)\\
15 & 0.16 & 0.15($-0.00$)\\
20 & 0.15 & 0.17($+0.02$)\\
25 & 0.16 & 0.21($+0.05$)\\
30 & 0.16 & 0.24($+0.08$)\\
  \bottomrule
 \end{tabular}
 \label{tab:hs_const_ann}
 \end{subtable}
 \hfill
 \begin{subtable}{.45\linewidth}
 \centering
 \captionsetup{font=footnotesize}
 \caption{Dynamic-load Test}
 \begin{tabular}{ccccc}
  \toprule
  $N_r$ & MLP + DPS & RKNN + DPS\\
  \midrule
1 & 0.15 & 0.16($+0.01$)\\
2 & 0.20 & 0.15($-0.05$)\\
3 & 0.15 & 0.12($-0.03$)\\
4 & 0.18 & 0.16($-0.02$)\\
5 & 0.18 & 0.14($-0.04$)\\
10 & 0.19 & 0.16($-0.03$)\\
15 & 0.17 & 0.16($-0.01$)\\
20 & 0.17 & 0.18($+0.02$)\\
25 & 0.16 & 0.19($+0.04$)\\
30 & 0.17 & 0.25($+0.08$)\\
  \bottomrule
 \end{tabular}
 \label{tab:hs_varnt_ann}
 \end{subtable}
\caption{MLP versus RKNN: Heat sink model, Relative error(\%)}
\label{tab:hs_mlpvsrknn}
\end{table}

\begin{table}[htbp]
\small
 \begin{subtable}{.45\linewidth}
 \centering
 \captionsetup{font=footnotesize}
 \caption{Constant-load Test}
 \begin{tabular}{ccccc}
  \toprule
  $N_r$ & MLP + DPS & RKNN + DPS\\
  \midrule
1 & 0.59 & 0.47($-0.13$)\\
2 & 0.27 & 0.24($-0.03$)\\
3 & 0.28 & 0.21($-0.07$)\\
4 & 0.27 & 0.19($-0.08$)\\
5 & 0.18 & 0.21($+0.02$)\\
10 & 0.14 & 0.19($+0.04$)\\
15 & 0.22 & 0.19($-0.03$)\\
20 & 0.13 & 0.19($+0.06$)\\
25 & 0.18 & 0.24($+0.06$)\\
30 & 0.16 & 0.21($+0.05$)\\
  \bottomrule
 \end{tabular}
 \label{tab:gr_const_ann}
 \end{subtable}
 \hfill
 \begin{subtable}{.45\linewidth}
 \centering
 \captionsetup{font=footnotesize}
 \caption{Dynamic-load Test}
 \begin{tabular}{ccccc}
  \toprule
  $N_r$ & MLP + DPS & RKNN + DPS\\
  \midrule
1 & 0.52 & 0.38($-0.13$)\\
2 & 0.19 & 0.25($+0.06$)\\
3 & 0.26 & 0.26($+0.01$)\\
4 & 0.15 & 0.25($+0.10$)\\
5 & 0.19 & 0.16($-0.03$)\\
10 & 0.19 & 0.21($+0.02$)\\
15 & 0.16 & 0.16($-0.01$)\\
20 & 0.13 & 0.19($+0.06$)\\
25 & 0.13 & 0.19($+0.06$)\\
30 & 0.12 & 0.20($+0.08$)\\
  \bottomrule
 \end{tabular}
 \label{tab:gr_const_ann}
 \end{subtable}
\caption{MLP versus RKNN: Gap radiation model, Relative error(\%)}
\label{tab:gr_mlpvsrknn}
\end{table}

\begin{table}[htbp]
\small
 \begin{subtable}{.45\linewidth}
 \centering
 \captionsetup{font=footnotesize}
 \caption{Constant-load Test}
 \begin{tabular}{ccccc}
  \toprule
  $N_r$ & MLP + DPS & RKNN + DPS\\
  \midrule
1 & 0.27 & 0.27($+0.01$)\\
2 & 0.27 & 0.13($-0.13$)\\
3 & 0.25 & 0.15($-0.10$)\\
4 & 0.24 & 0.20($-0.04$)\\
5 & 0.16 & 0.13($-0.03$)\\
10 & 0.11 & 0.18($+0.07$)\\
15 & 0.16 & 0.20($+0.05$)\\
20 & 0.12 & 0.16($+0.04$)\\
25 & 0.15 & 0.17($+0.02$)\\
30 & 0.13 & 0.17($+0.04$)\\
  \bottomrule
 \end{tabular}
 \label{tab:he_const_ann}
 \end{subtable}
 \hfill
 \begin{subtable}{.45\linewidth}
 \centering
 \captionsetup{font=footnotesize}
 \caption{Dynamic-load Test}
 \begin{tabular}{ccccc}
  \toprule
  $N_r$ & MLP + DPS & RKNN + DPS\\
  \midrule
1 & 0.22 & 0.18($-0.04$)\\
2 & 0.17 & 0.13($-0.04$)\\
3 & 0.27 & 0.15($-0.11$)\\
4 & 0.21 & 0.19($-0.02$)\\
5 & 0.15 & 0.15($+0.00$)\\
10 & 0.13 & 0.19($+0.06$)\\
15 & 0.17 & 0.21($+0.04$)\\
20 & 0.14 & 0.17($+0.03$)\\
25 & 0.15 & 0.18($+0.03$)\\
30 & 0.14 & 0.18($+0.04$)\\
  \bottomrule
 \end{tabular}
 \label{tab:he_const_ann}
 \end{subtable}
\caption{MLP versus RKNN: Heat exchanger model, Relative error(\%)}
\label{tab:he_mlpvsrknn}
\end{table}

In Table \ref{tab:gr_mlpvsrknn} and Table \ref{tab:he_mlpvsrknn}, we can observe that RKNN can predict the solution of reduced ODE system with acceptable accuracy. But there is no indication that RKNN predicts the system behavior better in any certain type of tests. In contrast, MLP is evaluated to be better at approximating larger-size ROM. One possible explanation is that compared to a simple MLP, an RKNN can be considered as a deeper network which is harder to train \cite{du2019gradient}. A solution could be to train large amount of networks for each architecture but with different randomized initialization. Then evaluating the average performance of all trained networks of each architecture. Here in this paper, due to the limitation of time and computational resource, we only train 10 networks for each type of architecture, which might not be sufficient to eliminate the random effect.

Although RKNN does present consistently smaller errors in our tests, it still has some advantages over the simple use of MLP. RKNN can be use to integrate a learned systems using a different time step in evaluation phase. According to Equation \ref{eq:rknn}, the knowledge learned by RKNN is the R.H.S. function $\boldsymbol{f}(\boldsymbol{y}_r;\boldsymbol{\mu})$ instead of the relation between $\boldsymbol{y}_r(t_j)$ and $\boldsymbol{y}_r(t_{j+1})$, we can use different $\tau$ to calculate the increment $\Delta \boldsymbol{y}_r(t_{j}) = \tau \boldsymbol{f}(\boldsymbol{y}_r(t_{j}); \boldsymbol{\mu}(t_j))$ without additional training process. Figure \ref{fig:hs_const_time_interval} - \ref{fig:he_varnt_time_interval} show how using different time interval can influence the prediction by both neural networks. An solution allowing MLP to predict with different step size is to include $\tau$ as a feature in the input of the MLP, according to Equation \ref{eq:direct_mlp} this can be written as:\\
\begin{equation}
    \boldsymbol{y}_{r,i} (t_{j+1}) \approx \boldsymbol{g}_{\text{dMLP}} (\boldsymbol{y}_{r,i} (t_{j}), \tau; \boldsymbol{\mu}(t_{j}))
\label{eq:direct_mlp_dt}
\end{equation}

However, in Equation \ref{eq:direct_mlp_dt}, $\tau$ contributes nonlinearly to the update of the state because of the nonlinear activation function in MLP, which is inconsistent to the pre-knowledge that the step size should contribute to the the update linearly $\tau \boldsymbol{f}(\boldsymbol{y}; \boldsymbol{\mu})$. Therefore, MLP can only approximate the relation, while RKNN exploits this pre-knowledge and embeds it into its structure.

\section{Conclusions}
\label{sec:conclusions}
In this paper, we proposed two modification for Model Order Reduction (MOR) based on Artificial Neural Network (ANN). The first modification is including more information of system parameters into training data by using Dynamic Parameter Sampling (DPS). The second is replacing conventional Multilayer Perceptron (MLP) with Runge-Kutta Neural Network (RKNN). The workflow of this technique is similar to conventional ANN-based MOR. Firstly we use POD reduction techniques to obtain an optimal reduced space from the snapshots and then, with the help of a neural network-based model, describes the dynamics of the system in such subspace. This approach is tested with two different data sampling strategies and neural network architectures respectively. The test cases are designed and simulated within NX 12.

As discussed in Section \ref{sec:results}, DPS shows positive influence on the quality of Reduced Order Model (ROM). Based on observation, adding DPS snapshots into dataset can enrich the dynamics contained in snapshots. This can help with constructing better POD basis. Moreover, neural network trained by DPS dataset is more sensitive to system's parameters. However, it is also seen that in some specific tests, network trained by DPS data performs worse than network trained by Static Parameter Sampling (SPS) data, which are usually found to be in constant-load tests.

As for applying RKNN instead of MLP, through error analysis, it is known that MLP and RKNN has similar capability to learning dataset in reduced space. However, some great potential of using RKNN to predict ODE system should be considered. For example, by learning the right hand side of the ODE system, RKNN becomes better architecture at enabling simulation with different step sizes from the training step size. This feature has been proven by results in Figure \ref{fig:hs_const_time_interval} - \ref{fig:he_varnt_time_interval}. Moreover, it is also mentioned in \cite{wang1998runge}, since the RKNN models the right hand side of ODE in its sub-networks directly, some known continuous relationship (physical laws) of the identiﬁed system can be incorporated into the RKNN to cut down the network size, speed up its learning, and enhance its long-term prediction capability.

Although this MOR framework was tested only on thermal-flow models, it could be easily adapted to other kinds of models, since it is non-intrusive and only relies on measured data to learn the dynamics of the system. However, there is still more work to be done before these methods can be applied to big scale and more general models.

Some further potential improvements for the approach discussed in this work include:
\begin{itemize}
    \item Investigating the influence of using SPS-DPS-mixed snapshots.
    \item Investigating different methods of constructing reduced space other than POD, e.g. Auto-encoder.
    \item Implementing implicit integration scheme with neural networks as in \cite{anderson1996comparison} allowing more accurate long-term prediction.
    \item Implementing integration scheme with multiple history states as part of input, which might also enable more stable long-term prediction.
\end{itemize}

\paragraph{Author Contributions}
Conceptualization: D.H., H.-J.B.; Methodology: D.H., Q.Z. and J.L.; Data curation: Q.Z.; Data visualisation: Q.Z.; Writing original draft: D.H., Q.Z. and J.L.; All authors approved the final submitted draft.

\bibliographystyle{apalike}
\bibliography{references}

\newpage
\appendix
\section{Additional Figures}
\label{appendix}

\begin{figure}[htbp]
  \centering
  \begin{subfigure}[b]{0.45\textwidth}
  \centering
  \includegraphics[width=\textwidth]{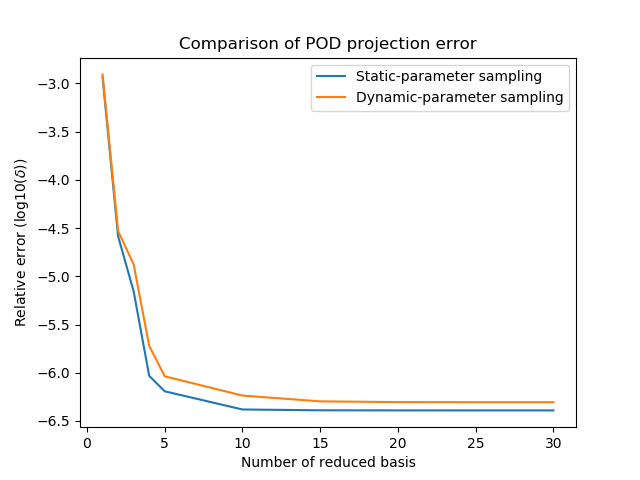}
  \captionsetup{font=footnotesize}
  \caption{Constant load test}
  \label{fig:poder_heatsink_const}
  \end{subfigure}
  \hfill
  \begin{subfigure}[b]{0.45\textwidth}
  \centering
  \includegraphics[width=\textwidth]{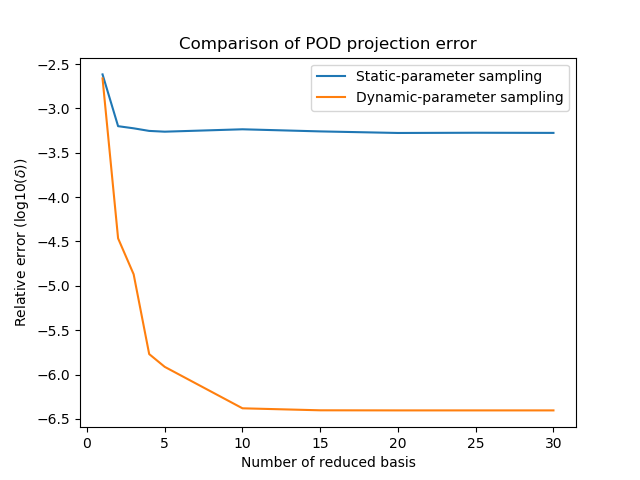}
  \captionsetup{font=footnotesize}
  \caption{Dynamic load test}
  \label{fig:poder_heatsink_varnt}
  \end{subfigure}
  \hfill
\caption{Projection error: heat sink model}
\label{fig:pod_error_heatsink}
\end{figure}

\begin{figure}[htbp]
  \centering
  \begin{subfigure}[b]{0.45\textwidth}
  \centering
  \includegraphics[width=\textwidth]{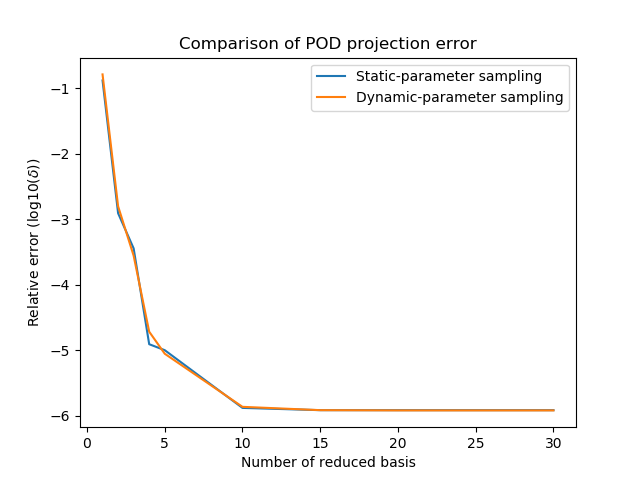}
  \captionsetup{font=footnotesize}
  \caption{Constant load test}
  \label{fig:poderr_radiation_const}
  \end{subfigure}
  \hfill
  \begin{subfigure}[b]{0.45\textwidth}
  \centering
  \includegraphics[width=\textwidth]{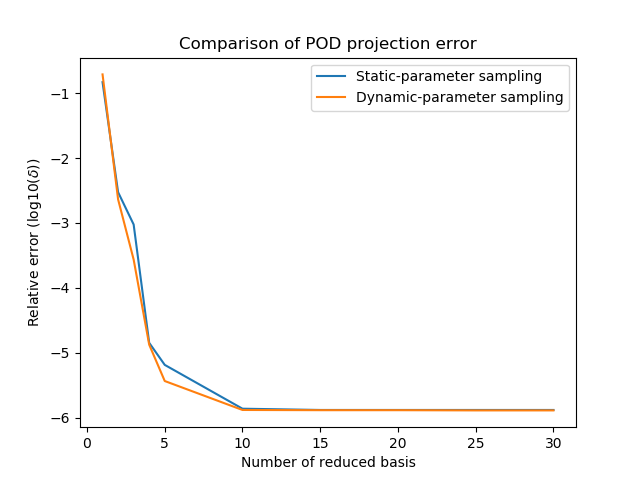}
  \captionsetup{font=footnotesize}
  \caption{Dynamic load test}
  \label{fig:poderr_radiation_varnt}
  \end{subfigure}
  \hfill
\caption{Projection error: gap radiation model}
\label{fig:pod_error_radiation}
\end{figure}

\begin{figure}[htbp]
  \centering
  \begin{subfigure}[b]{0.45\textwidth}
  \centering
  \includegraphics[width=\textwidth]{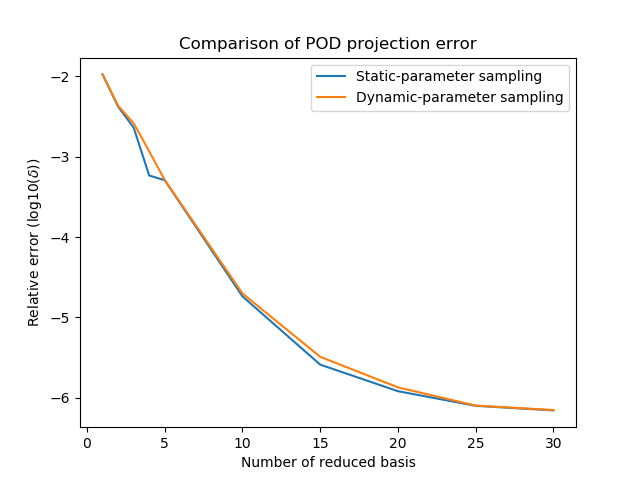}
  \captionsetup{font=footnotesize}
  \caption{Constant load test}
  \label{fig:poderr_radiation_const}
  \end{subfigure}
  \hfill
  \begin{subfigure}[b]{0.45\textwidth}
  \centering
  \includegraphics[width=\textwidth]{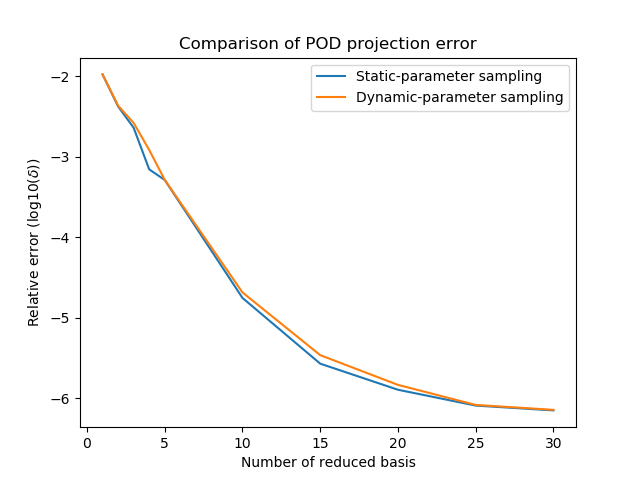}
  \captionsetup{font=footnotesize}
  \caption{Dynamic load test}
  \label{fig:poderr_exchanger_varnt}
  \end{subfigure}
  \hfill
\caption{Projection error: heat exchanger model}
\label{fig:pod_error_radiation}
\end{figure}

\begin{figure}[htbp]
  \centering
  \begin{subfigure}[b]{0.45\textwidth}
  \centering
  \includegraphics[width=\textwidth]{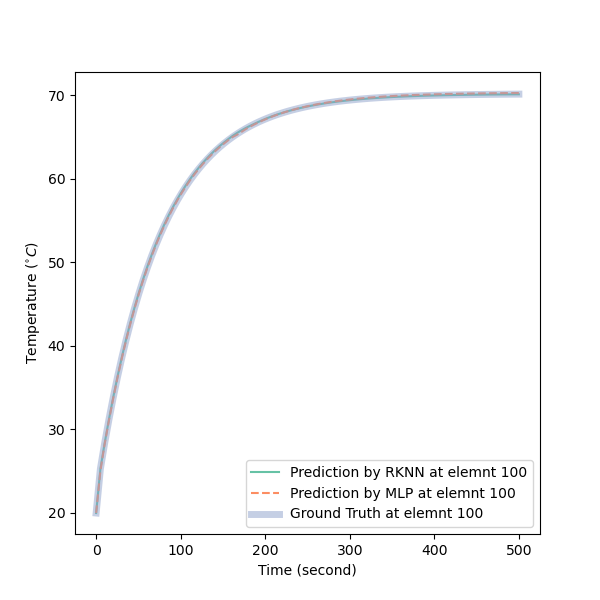}
  \captionsetup{font=footnotesize}
  \caption{Time interval: dt=5s}
  \label{fig:hs_const_dt}
  \end{subfigure}
  \hfill
  \begin{subfigure}[b]{0.45\textwidth}
  \centering
  \includegraphics[width=\textwidth]{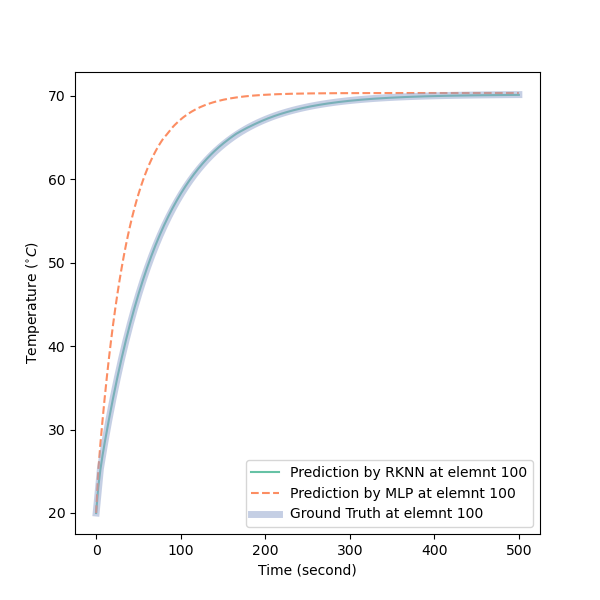}
  \captionsetup{font=footnotesize}
  \caption{Time interval: dt=2.5s}
  \label{fig:hs_const_0.5dt}
  \end{subfigure}
  \hfill
\caption{ANN predictions: constant-load test, heat sink model}
\label{fig:hs_const_time_interval}
\end{figure}

\begin{figure}[htbp]
  \centering
  \begin{subfigure}[b]{0.45\textwidth}
  \centering
  \includegraphics[width=\textwidth]{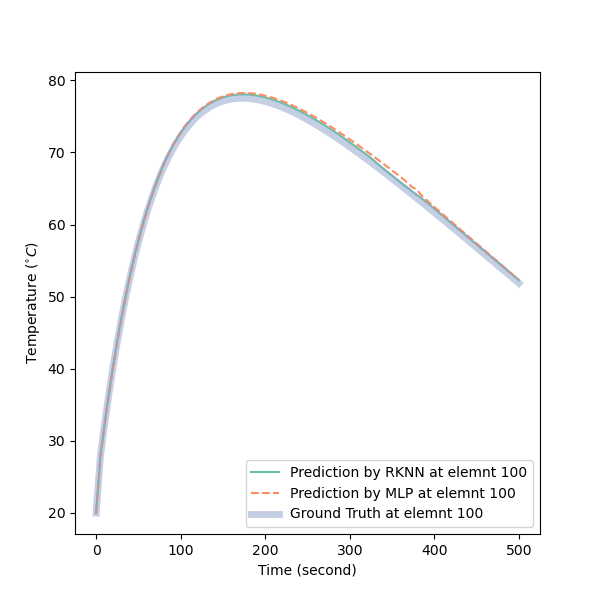}
  \captionsetup{font=footnotesize}
  \caption{Time interval: dt=5s}
  \label{fig:hs_varnt_dt}
  \end{subfigure}
  \hfill
  \begin{subfigure}[b]{0.45\textwidth}
  \centering
  \includegraphics[width=\textwidth]{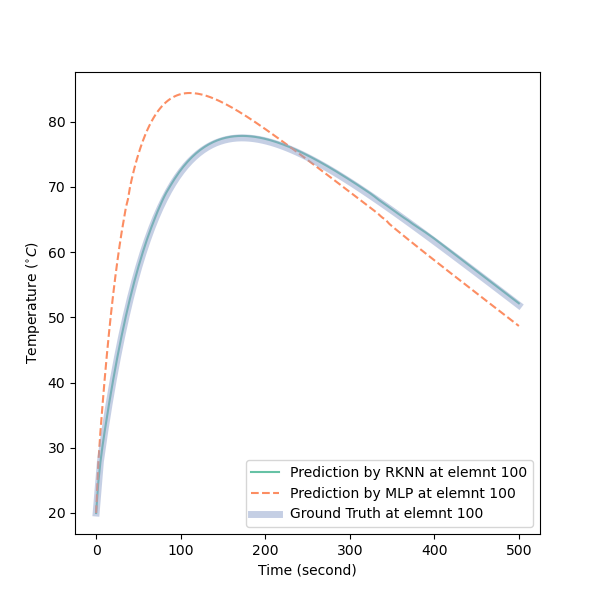}
  \captionsetup{font=footnotesize}
  \caption{Time interval: dt=2.5s}
  \label{fig:hs_varnt_0.5dt}
  \end{subfigure}
  \hfill
\caption{ANN predictions: dynamic-load test, heat sink model}
\label{fig:hs_varnt_time_interval}
\end{figure}

\begin{figure}[htbp]
  \centering
  \begin{subfigure}[b]{0.45\textwidth}
  \centering
  \includegraphics[width=\textwidth]{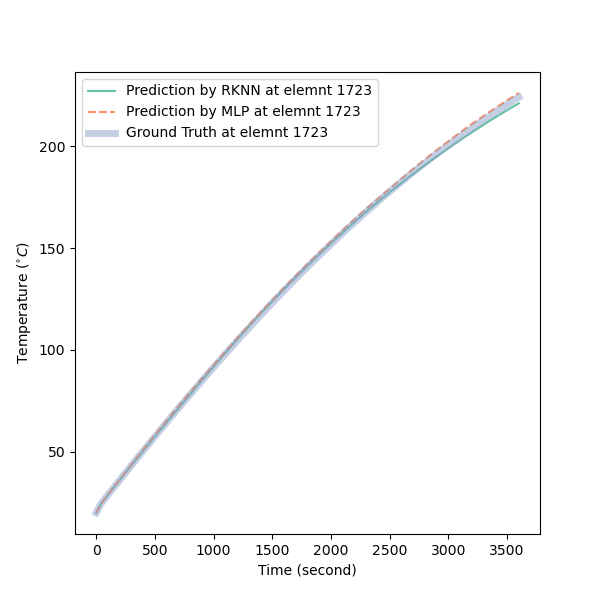}
  \captionsetup{font=footnotesize}
  \caption{Time interval: dt=36s}
  \label{fig:gr_const_dt}
  \end{subfigure}
  \hfill
  \begin{subfigure}[b]{0.45\textwidth}
  \centering
  \includegraphics[width=\textwidth]{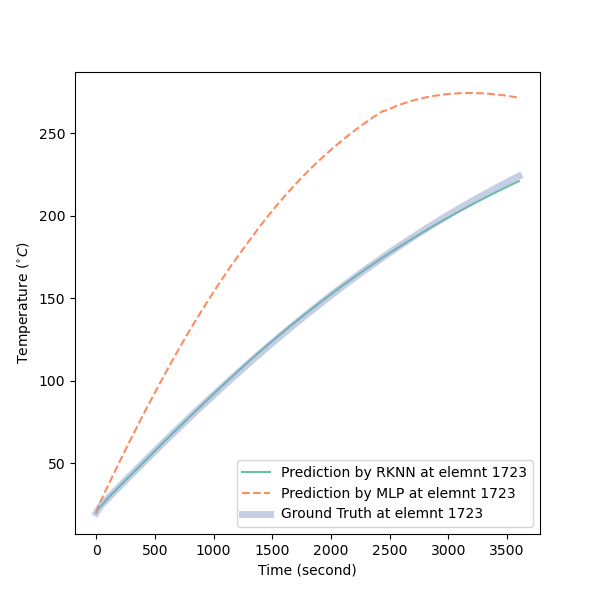}
  \captionsetup{font=footnotesize}
  \caption{Time interval: dt=18s}
  \label{fig:gr_const_0.5dt}
  \end{subfigure}
  \hfill
\caption{ANN predictions: constant-load test, gap radiation model}
\label{fig:gr_const_time_interval}
\end{figure}

\begin{figure}[htbp]
  \centering
  \begin{subfigure}[b]{0.45\textwidth}
  \centering
  \includegraphics[width=\textwidth]{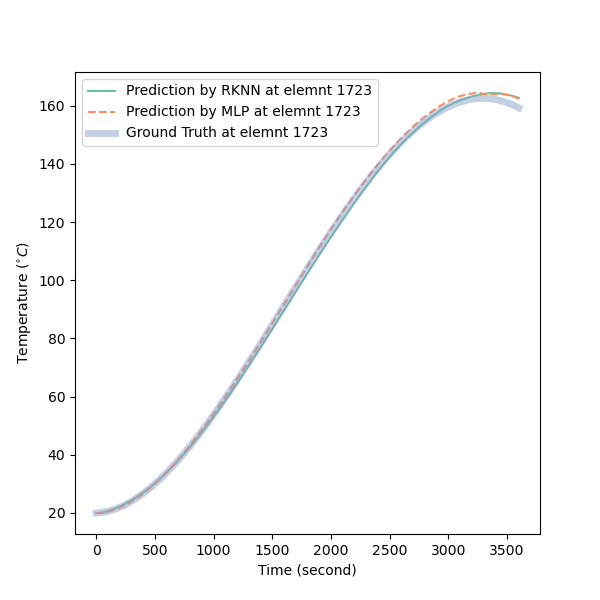}
  \captionsetup{font=footnotesize}
  \caption{Time interval: dt=36s}
  \label{fig:gr_varnt_dt}
  \end{subfigure}
  \hfill
  \begin{subfigure}[b]{0.45\textwidth}
  \centering
  \includegraphics[width=\textwidth]{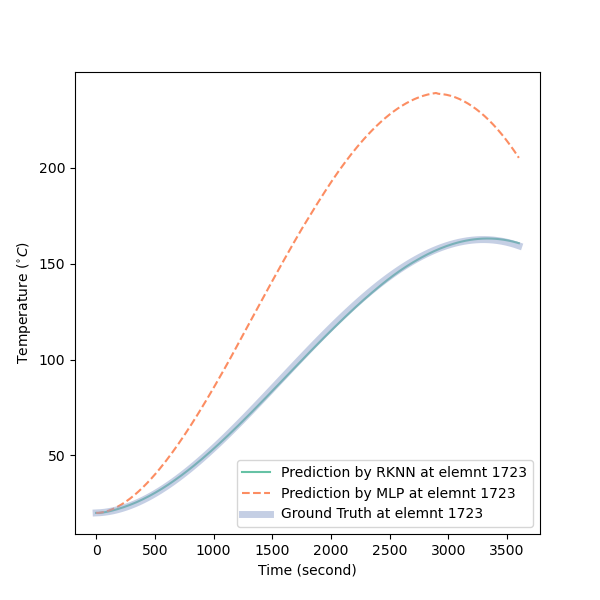}
  \captionsetup{font=footnotesize}
  \caption{Time interval: dt=18s}
  \label{fig:gr_varnt_0.5dt}
  \end{subfigure}
  \hfill
\caption{ANN predictions: dynamic-load test, gap radiation model}
\label{fig:gr_varnt_time_interval}
\end{figure}

\begin{figure}[htbp]
  \centering
  \begin{subfigure}[b]{0.45\textwidth}
  \centering
  \includegraphics[width=\textwidth]{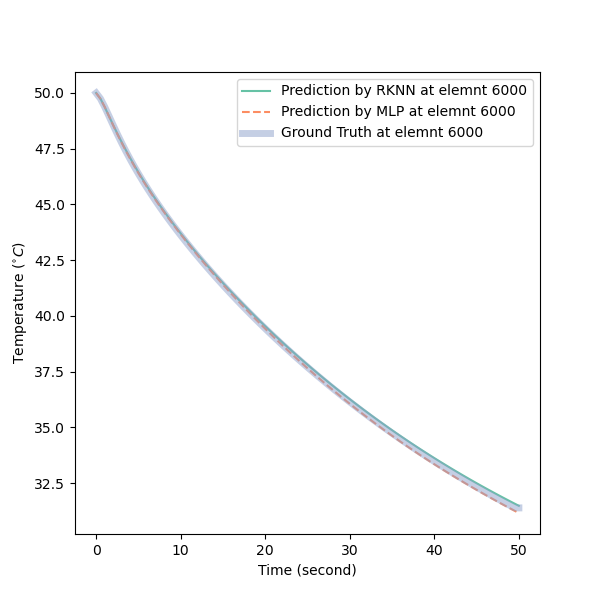}
  \captionsetup{font=footnotesize}
  \caption{Time interval: dt=0.5s}
  \label{fig:he_const_dt}
  \end{subfigure}
  \hfill
  \begin{subfigure}[b]{0.45\textwidth}
  \centering
  \includegraphics[width=\textwidth]{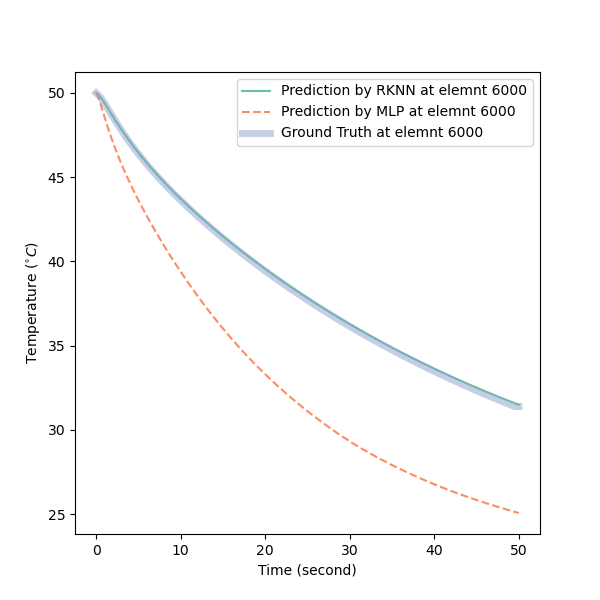}
  \captionsetup{font=footnotesize}
  \caption{Time interval: dt=0.25s}
  \label{fig:he_const_0.5dt}
  \end{subfigure}
  \hfill
\caption{ANN predictions: constant-load test, heat exchanger model}
\label{fig:he_const_time_interval}
\end{figure}

\begin{figure}[htbp]
  \centering
  \begin{subfigure}[b]{0.45\textwidth}
  \centering
  \includegraphics[width=\textwidth]{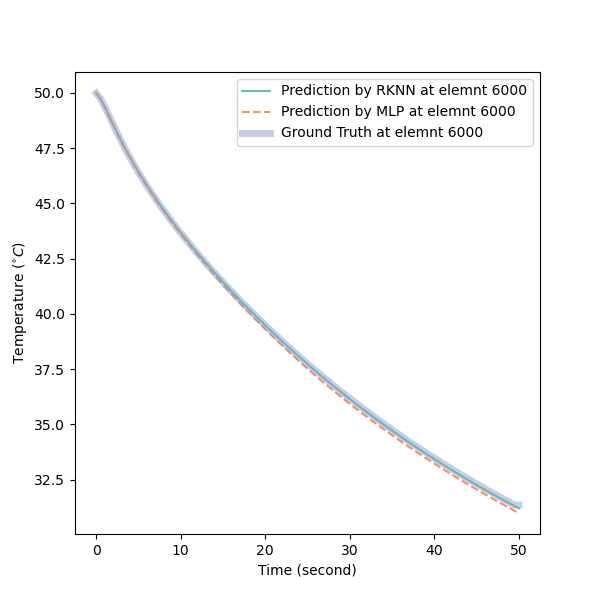}
  \captionsetup{font=footnotesize}
  \caption{Time interval: dt=0.5s}
  \label{fig:he_varnt_dt}
  \end{subfigure}
  \hfill
  \begin{subfigure}[b]{0.45\textwidth}
  \centering
  \includegraphics[width=\textwidth]{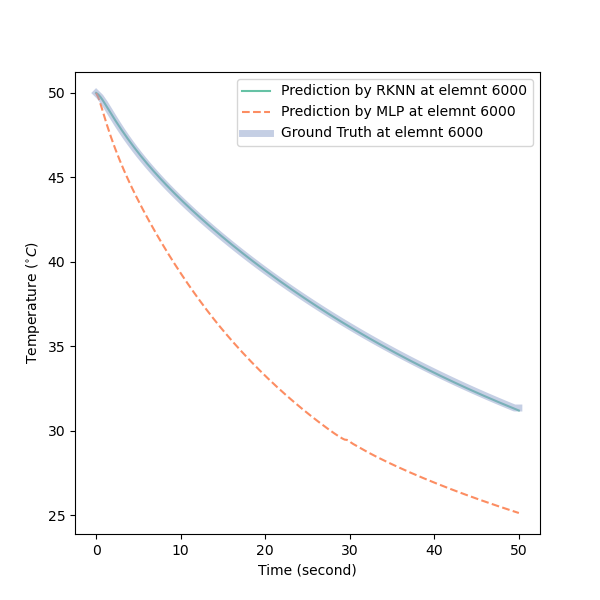}
  \captionsetup{font=footnotesize}
  \caption{Time interval: dt=0.25s}
  \label{fig:he_varnt_0.5dt}
  \end{subfigure}
  \hfill
\caption{ANN predictions: dynamic-load test, heat exchanger model}
\label{fig:he_varnt_time_interval}
\end{figure}

\end{document}